\documentclass[12pt]{article}

\pdfoutput=1

\usepackage{amssymb,amsmath,epsfig,amscd,eucal,psfrag}

\newtheorem{thm}{Theorem}[section]
\newtheorem{defn}[thm]{Definition}
\newtheorem{lem}[thm]{Lemma}
\newtheorem{prop}[thm]{Proposition}

\newtheorem{rem}[thm]{Remark}
\newtheorem{cor}[thm]{Corollary}

\newtheorem{letterthm}{Theorem}

\newenvironment{pf}{\par\medskip\noindent{\em Proof. }}{\hfill $\square$\par\medskip}
\newenvironment{pfof}[1]{\par\medskip\noindent{\em Proof of #1. }}{\hfill $\square$\par\medskip}

\newcommand{\F}{\mathbb{F}}

\newcommand{\R}{\mathbb{R}}
\newcommand{\Z}{\mathbb{Z}}

\newcommand{\ICE}{\mathrm{ICE}}
\newcommand{\Axis}{\mathrm{Axis}}

\newcommand{\Hom}{\mathrm{Hom}}

\newcommand{\curlyL}{\mathcal{L}}

\title{Hall's Theorem for limit groups}
\author{Henry Wilton\footnote{Partially supported by an EPSRC student scholarship and by a post-doctoral fellowship at the Hebrew University of Jerusalem, Israel.}}

\date{2nd May 2007}

\begin{document}
\maketitle

\begin{abstract}
A celebrated theorem of Marshall Hall Jr.\ implies that finitely generated free groups are subgroup separable and that all of their finitely generated subgroups are retracts of finite-index subgroups.  We use topological techniques inspired by the work of Stallings to prove that all limit groups share these two properties.  This answers a question of Sela.
\end{abstract}

Limit groups are finitely presented groups that arise naturally in many different aspects of the study of finitely generated (non-abelian) free groups.  Perhaps their most satisfying characterization is as the closure of the set of free groups in the topology on marked groups that arose from the work of M.~Gromov and R.~Grigorchuk (see \cite{CG04}).  Limit groups admit a hierarchical decomposition in which the basic building blocks are free groups, free abelian groups and the fundamental groups of surfaces of Euler characteristic less than -1.  Therefore, it seems natural to try to generalize properties of these ubiquitous classes of groups to limit groups.

Much of the recent work on limit groups has been motivated by the fundamental role that they play in the study of $\Hom(G,\F)$, the variety of homomorphisms from a finitely generated group $G$ to a free group $\F$, and in the first-order logic of the free group. One can associate to a group $G$ the \emph{elementary theory} of $G$, the set of sentences in first-order logic that hold in $G$. The elementary theory contains the \emph{existential theory}, which consists of those sentences that use only one, existential, quantifier.  From this point of view limit groups are precisely the finitely generated groups with the same existential theory as a free group \cite{Rem89}.  In \cite{Se1,Se2,Se3,Se4,Se5a,Se5b,Se6}, Zlil Sela solved a famous problem of Alfred Tarski by classifying the \emph{elementarily free groups}, those limit groups with the same elementary theory as a free group.  O.~Kharlampovich and A.~Miasnikov also announced a solution to Tarski's Problem in \cite{KM98a,KM98b,KM3,KM4,KM5}.

Work on Tarski's Problem has led to the development of a powerful structure theory for limit groups.  One form of this is given by a theorem of Kharlampovich and Miasnikov \cite{KM98b} (see also \cite{CG04}).  Let $\ICE$ be the smallest set of groups containing all finitely generated free groups that is closed under extending centralizers.  Limit groups are precisely the finitely generated subgroups of groups in $\ICE$ . We make extensive use of this characterization.

Marshall Hall Jr.\ proved in \cite{Hall49} that every finitely generated subgroup of a finitely generated free group is a free factor in a finite-index subgroup, and that this finite-index subgroup can be chosen to exclude any element not in the original subgroup.  As a consequence, finitely generated free groups are \emph{subgroup separable} (also known as \emph{LERF}); that is, every finitely generated subgroup is closed in the profinite topology.   This is a strong algebraic condition that implies, for instance, that the generalized word problem is solvable.  In the 1970s it became apparent that subgroup separability has a natural topological interpretation.  Peter Scott used hyperbolic geometry in \cite{Sc1} to prove that surface groups are subgroup separable, while J.~R.~Stallings exploited the topology of graphs to reprove Hall's Theorem and other properties of free groups in \cite{St}.  Stallings' techniques have been extended by, among others, Rita Gitik \cite{Gi1}, who proved that, in certain circumstances, the amalgamated product of a subgroup separable group with a free group over a cyclic subgroup is subgroup separable, and D.~T.~Wise \cite{Wise1}, who classified the subgroup separable graphs of free groups with cyclic edge groups. Sela \cite{SeQ} asked if limit groups are subgroup separable.  In \cite{Wilt1}, the author answered Sela's question in the affirmative for elementarily free groups. Here we extend that result.

\begin{letterthm}[Corollary \ref{LGs are subgroup separable}]
\label{Theorem A} Limit groups  are subgroup separable.
\end{letterthm}

It follows that limit groups have solvable generalized word problem.  We also prove the stronger theorem that $\ICE$ groups are coset separable with respect to vertex groups.

Another strong consequence of Hall's Theorem is that every finitely generated subgroup $H$ of a free group $F$ is a \emph{virtual retract}; that is, $H$ is a retract of some finite-index subgroup of $F$. This property is referred to as \emph{local retractions} in \cite{LR2}, in which connections between virtual retracts and 3-manifold topology are explored. Combining Theorem \ref{Theorem A} with Theorem 3.1 of \cite{BrHo} it follows that every cyclic subgroup of a limit group is a virtual retract.  We show here that limit groups virtually retract onto all their finitely generated subgroups.

\begin{letterthm}[Corollary \ref{LGs have local retractions}]
\label{Theorem B} Let $G$ be a limit group and $H$ a finitely generated subgroup of $G$.  Then there exists a finite-index subgroup $K$ of $G$ containing $H$ and a retraction $K \to H$.
\end{letterthm}

As a corollary one obtains a purely topological proof that finitely generated subgroups of limit groups are quasi-isometrically embedded.  See subsection \ref{QI subgroups} for more details.

The paper is organized as follows.  In section 1 the notions of subgroup separability and local retractions are introduced.  An example is given of how Stallings' techniques can be used to prove a strong theorem about finitely generated free groups.  In order to extend these ideas, we need graphs of spaces in the spirit of \cite{SW}. We briefly outline some salient features of the theory of limit groups and their structure. In section 2 we develop the language of elevations (used by Wise in \cite{Wise1}) and pre-covers (analogous to those used by Gitik in \cite{Gi1}). Using these we can state a precise generalization of Stallings' ideas to the context of graphs of spaces (Theorem \ref{Stallings}), which we call Stallings' Principle.  There are still substantial obstacles to applying this to limit groups, and section 3 is devoted to overcoming them.  Specifically, we define tame coverings, in which we have strong control over elevations of loops.  The main technical result, Theorem \ref{Main theorem}, asserts that coverings of $K(G,1)$s for groups in $\ICE$ are tame over certain sets of loops. Theorems \ref{Theorem A} and \ref{Theorem B} follow quickly.

This paper makes use of the results about pre-covers and elevations in \cite{Wilt1}, whither we often refer the reader for proofs.  However, it should be emphasized that Theorems \ref{Theorem A} and \ref{Theorem B} are independent of the main theorem of \cite{Wilt1} (Theorem 5.13).  Theorem \ref{Theorem A} was announced in \cite{Wilt3}, with a citation of an earlier version of this paper entitled `Limit Groups are Subgroup Separable'.

Throughout, if $g$ and $h$ are group elements then $g^h$ denotes $h^{-1}gh$.

\subsection*{Acknowledgements}

I would like to thank William Dison, Andrew Harkins and Michael Tweedale for patiently listening to rehearsals of the ideas of this work, Alan Reid for drawing local retractions to my attention, and particularly Martin Bridson for his many insightful comments and suggestions.

\section{Fundamental notions}

\subsection{Separability of subgroups and double cosets}

\begin{defn}
A subgroup $H\subset G$ is \emph{separable} if it is an intersection of finite-index subgroups of $G$; equivalently, for every $g\in G\smallsetminus H$ there exists a finite-index subgroup $K\subset G$ such that $H\subset K$ but $g\notin K$.

A group is \emph{subgroup separable} if every finitely generated subgroup is separable.
\end{defn}
Note that if $G$ is subgroup separable and $H\subset G$ is a subgroup then $H$ is also subgroup separable.

A group $G$ has solvable \emph{generalized word problem} if there exists an algorithm that, given an element $g\in G$ and a finite subset $S\subset G$, determines whether or not $g\in\langle S\rangle$.   Just as residually finite groups have solvable word problem, so subgroup separable groups have solvable generalized word problem.

We will also consider the separability of double cosets.

\begin{defn}
If $H$ and $H'$ are subgroups of $G$ and $g\in G$, the double coset $HgH'$ is separable if, whenever $HgH'\neq HhH'$, there exists a finite-index subgroup $K'\subset G$, containing $H'$, such that $HgK'\neq HhK'$.

The group $G$ is called \emph{coset separable with respect to $H$} if, for every finitely generated subgroup $H'\subset G$ and every $g\in G$, the double coset $HgH'$ is separable.
\end{defn}

Subgroup separability has a beautiful topological reformulation.  We work in the category of combinatorial complexes and combinatorial maps.  If $G$ is the fundamental group of a complex $X$ then the separability of the subgroup $H$ can be reformulated in terms of a condition on the covering of $X$ corresponding to $H$.

\begin{lem}[Scott's Criterion]\label{Topological interpretation}
Let $X$ be a connected complex, $H$ a subgroup of $\pi_1(X)$ and $p:X^H\to X$ the (connected) covering corresponding to $H$.  The following are equivalent.
\begin{enumerate}
\item The subgroup $H$ is separable.
\item For any finite subcomplex $\Delta\subset X^H$ there exists a finite-sheeted intermediate covering
$$
X^H\to \hat{X}\to X
$$
so that $\Delta$ embeds in $\hat{X}$.
\end{enumerate}
\end{lem}
\begin{pf}
Fix basepoints $x_0\in X$ and $x'_0\in p^{-1}(x_0)$.  Suppose $H$ is separable and $\Delta\subset X^H$ is a finite subcomplex.  For each cell $\sigma\subset\Delta$ fix a choice of barycentre $a_\sigma\in \sigma$, in such a way that if cells $\sigma$ and $\sigma'$ have the same image in $X$ then $p(a_\sigma)=p(a_{\sigma'})$.  Choose paths $\alpha_\sigma$ in $X^H$ from $x'_0$ to $a_\sigma$ and denote by $\bar{\alpha}_\sigma$ the reverse path from $a_\sigma$ to $x'_0$.  For each distinct pair of cells $\sigma,\sigma'\subset \Delta$ with $p(\sigma)=p(\sigma')$ the concatenation of paths $\gamma_{\sigma,\sigma'}=p\circ\alpha_{\sigma}\cdot p\circ\bar{\alpha}_{\sigma'}$ lies in $\pi_1(X,x_0)$ but not in $H$.  Since $H$ is separable, for each such pair $\sigma,\sigma'$ there exists a finite-index subgroup $K_{\sigma,\sigma'}\subset\pi_1(X)$ with $H\subset K_{\sigma,\sigma'}$ but $g_{\sigma,\sigma'}\notin K_{\sigma,\sigma'}$.  Define a finite-sheeted covering $\hat{X}\to X$ by
\[
\pi_1(\hat{X})=\bigcap_{\sigma,\sigma'} K_{\sigma,\sigma'}
\]
where the intersection ranges over each all pairs of cells $\sigma,\sigma'\subset\Delta$ with $p(\sigma)=p(\sigma')$.  Then $\Delta$ embeds in $\hat{X}$ by construction.

Conversely, suppose every finite subcomplex of $X^H$ embeds in a finite-sheeted cover and let $g\in \pi_1(X)\smallsetminus H$.  Fix a combinatorial loop representing $g$ and consider its based lift $g':[0,1]\to X^H$ with $g'(0)=x'_0$.  Since $g\notin H$, $g'(1)\neq x'_0$.  By the supposition, there exists a finite-sheeted intermediate covering $\hat{X}$ such that the image of $g'$ embeds in $\hat{X}$.  In particular, $g\notin\pi_1(\hat{X})$.  Hence, $H$ is separable.
\end{pf}

\subsection{Local retractions}

A group $G$ \emph{retracts onto} a subgroup $H$ if the inclusion map $H\hookrightarrow G$ has a left-inverse $\rho$. In this case the map $\rho$ is called a \emph{retraction} and the subgroup $H$ is called a \emph{retract}.

\begin{defn}[\cite{LR2}]
A group $G$ \emph{virtually retracts onto} a subgroup $H$ if there exists a finite-index subgroup $K\subset G$ containing $H$ and a retraction $K\to H$.  If $G$ virtually retracts onto all its finitely generated subgroups then $G$ is said to have \emph{local retractions} or \emph{property LR}.
\end{defn}

Finitely generated abelian groups have local retractions.  Hall's Theorem asserts that every finitely generated subgroup $H$ of a free group $F$ is a free factor in a subgroup $K\subset F$ of finite index; in particular, $H$ is a retract of $K$, so $F$ has local retractions.  Notice that having local retractions passes to subgroups.

Having local retractions is a strong property, with many deep consequences. For more information, see \cite{LR2}.  We will mention just one such consequence.  Recall that a map of metric spaces $f:X\to Y$ is a \emph{quasi-isometric embedding} if there exist constants $\lambda$ and $\epsilon$ such that, for all $x_1,x_2\in X$,
$$
\lambda^{-1}d_X(x_1,x_2)-\epsilon\leq d_Y(f(x_1),f(x_2))\leq\lambda d_X(x_1,x_2)+\epsilon.
$$
A finitely generated subgroup $H$ of a finitely generated group $G$ is \emph{quasi-isometrically embedded} if, for some choice of finite generating sets and associated word metrics, the inclusion map $H\hookrightarrow G$ is a quasi-isometric embedding. This is in fact independent of the choices of generating sets.  Since the word metric on a group is only well-defined up to quasi-isometry, this is a natural notion of good geometric behaviour for a subgroup.

\begin{lem}\label{LR implies qi embedded subgroups}
If a finitely generated group $G$ virtually retracts onto a finitely generated subgroup $H$ then $H$ is quasi-isometrically embedded.
\end{lem}
\begin{pf}
There exists finite-index $K\subset G$ and a retraction $\rho:K\to H$. Since $K$ decomposes as $H\ker\rho$ there is a finite generating set for $K$ of the form $S\cup T$ where $H=\langle S\rangle$ and $T\subset\ker\rho$.  In the word metric with respect to this generating set, the inclusion $H\subset K$ is an isometric embedding. But $K\hookrightarrow G$ is a quasi-isometry, so $H\hookrightarrow G$ is a quasi-isometric embedding.
\end{pf}

\subsection{The motivating idea}\label{Motivating example subsection}

The proof of our main technical result, Theorem \ref{Main theorem}, is inspired by Stallings' famous proof of Hall's Theorem \cite{St} which we outline here.

\begin{thm}\label{Immersions of graphs}
Let $X'\to X$ be an immersion of finite graphs.  Then $X'$ embeds in some finite graph $\hat{X}$ and there exists a covering map $\hat{X}\to X$ extending $X'\to X$.
\end{thm}
We postpone the proof of this to subsection \ref{Subsection Stallings' Principle}.

\begin{cor}[\cite{Hall49}]\label{Hall's Theorem for free groups}
Finitely generated free groups are subgroup separable and have local retractions.
\end{cor}
\begin{pf}
Let $F$ be a finitely generated free group and $H$ a finitely generated subgroup.  Realize $F$ as the fundamental group of some finite graph $X$, in which case $H$ corresponds to a covering $X^H\to X$.  Let $\Delta$ be a finite subgraph of $X^H$.  Since $H$ is finitely generated, there exists a finite connected subgraph $X'\subset X^H$ so that the inclusion map is a $\pi_1$-isomorphism. Enlarging $X'$ if necessary, it can be assumed that $\Delta\subset X'$.  The restriction of the covering map $X'\to X$ is an immersion of finite graphs.  By Theorem \ref{Immersions of graphs}, $X'$ can be completed to a finite-sheeted cover $\hat{X}$ into which $\Delta$ embeds.  Since $X'$ is a subgraph of $\hat{X}$, $H$ is a free factor and, in particular, a retract of $\pi_1(\hat{X})$.
\end{pf}

We aim to prove a version of Theorem \ref{Immersions of graphs} in the context of graphs of spaces (or, equivalently, graphs of groups---see subsection \ref{Subsection graphs of spaces}). Three general theorems have been proved about subgroup separability of graphs of groups.  One concerns the amalgamation of a subgroup separable group and a free group \cite{Gi1}.  Another concerns graphs of free groups \cite{Wise1}.  In \cite{Wilt1}, we proved that elementarily free groups are subgroup separable.  (Elementarily free groups can be constructed inductively by gluing surfaces along their boundaries.)  In all these theorems, at least one vertex group of the graph of groups is free.  Free groups play such a key role because Theorem \ref{Immersions of graphs} provides great flexibility in constructing finite-index subgroups of free groups with special properties.

The flexibility required in \cite{Gi1}, \cite{Wilt1} and \cite{Wise1} is similar to the following corollary of Theorem \ref{Immersions of graphs}.  A collection $g_1,\ldots,g_n$ of elements of a group $G$ is \emph{independent} if, whenever there exists $h\in G$ such that $g_i^h$ and $g_j$ commute, then in fact $i=j$.

\begin{cor}\label{Motivating example}
Let $X$ be a graph and let $H\subset\pi_1(X)$ be finitely generated. Let $X^H\to X$ be the covering corresponding to $H$ and $\Delta\subset X^H$ a finite subcomplex. Let $\{\gamma_i\}$ be a finite independent set of elements of $\pi_1(X)$ that each generate a maximal cyclic subgroup. For each $\gamma_i$, consider a finite collection of conjugates $\{\gamma_i^{g_j}\}$, such that $\langle \gamma_i^{g_j} \rangle\cap H=1$ for each $i$ and $j$. Then for all sufficiently large positive integers $d$ there exists a finite-sheeted intermediate covering
$$
X^H\to\hat{X}\to X
$$
such that:
\begin{enumerate}
\item no distinct $\gamma_i^{g_j}$ and $\gamma_i^{g_{j'}}$ are conjugate in $\pi_1(\hat{X})$;
\item $\langle \gamma_i^{g_j} \rangle\cap \pi_1(\hat{X})=\langle (\gamma_i^{g_j})^d\rangle$ for each $i$ and $j$;
\item $\Delta$ embeds into $\hat{X}$;
\item there exists a retraction $\rho:\pi_1(\hat{X})\to H$ and $\rho(\gamma_i^{g_j})=1$ for each $i$ and $j$.
\end{enumerate}
\end{cor}
Note that Corollary \ref{Motivating example} implies Corollary \ref{Hall's Theorem for free groups}.

\begin{pf}
Fix a basepoint $x\in X$ and let $(\tilde{X},\tilde{x})\to (X,x)$ be the universal covering. Since $H$ is finitely generated there exists some connected finite subgraph $X'\subset X^H$ such that the inclusion map is a $\pi_1$-isomorphism.  Enlarging $X'$ we can assume that it contains $\Delta$ and the basepoint $x'\in X^H$.  Note that $X^H\smallsetminus X'$ is a forest.

Denote the universal covering map $\R\to S^1$ by $p$.  Identifying each loop $\gamma_i^{g_j}$ with a based map $(S^1,p(0))\to X$, each concatenation $\gamma_i^{g_j}\circ p$ lifts to a based map
$$
\gamma_j^H:(\R,0)\to (X^H,x')
$$
and further to
$$
\tilde{\gamma}_j:(\R,0)\to (\tilde{X},\tilde{x}),
$$
where we have dropped the $i$'s for notational convenience.  Since $\langle \gamma_i^{g_j} \rangle\cap H=1$ the $\gamma_j^H$ are proper.  Each $\tilde{\gamma}_j$ contains $\Axis(\gamma_i^{g_j})$ in its image .

If the intersection of the images of $\gamma^H_j$ and $\gamma^H_{j'}$ contains an infinite ray then, for some $h\in H$, the intersection
$$
h\Axis(\gamma_i^{g_j})\cap\Axis(\gamma_{i'}^{g_{j'}})
$$
contains an infinite ray.  Therefore, for some integers $m$ and $n$, the element
$$
h(\gamma_i^{g_j})^mh^{-1}(\gamma_{i'}^{g_{j'}})^n\in \pi_1(X)
$$
fixes a point in $\tilde{X}$ and so is trivial, whence $i=i'$ and $j=j'$.  Enlarging $X'$ still further we can also assume that, for all sufficiently large $d$:
\begin{enumerate}
\item each $\gamma^H_j$ restricts to an arc
$$
\gamma'_j:[a_j,b_j]\to X';
$$
of combinatorial length $d$ times the combinatorial length of $\gamma_i$;
\item for any distinct $\gamma_i^{g_j}$ and $\gamma_{i'}^{g_{j'}}$,
$$
|\{\gamma'_j(a_j),\gamma'_j(b_j)\}\cup\{\gamma'_{j'}(a_{j'}),\gamma'_{j'}(b_{j'})\}|=4.
$$
\end{enumerate}
By identifying each $\gamma'_j(a_j)$ with $\gamma'_j(b_j)$ we obtain an immersion $\bar{X}\to X$ such that $\langle\gamma_i^{g_j}\rangle\cap\pi_1(\bar{X})=\langle(\gamma_i^{g_j})^d\rangle$ for every $i$ and $j$.  Applying Theorem \ref{Immersions of graphs} to $\bar{X}$ gives the required cover $\hat{X}$. Note that $X'$ is a subgraph of $\hat{X}$ so $H$ is a free factor of $\pi_1(\hat{X})$ and, in particular, a retract.
\end{pf}

In the terminology of section \ref{Constructing finite-sheeted covers section}, Corollary \ref{Motivating example} asserts that the cover $X^H$ is tame over the set $\{\gamma_i\}$.  The lifts $\gamma^H_j$ are elevations and $X'$ is a pre-cover.

The proof of Corollary \ref{Motivating example} makes use of the fact that the action of $\pi_1(X)$ on the tree $\tilde{X}$ is free.  In extending these ideas to limit groups, the key observation is that a similar construction can be made when we make the weaker assumption that an action on a tree is \emph{acylindrical} (see Lemma \ref{Ensuring disparity}).

\subsection{Graphs of spaces}\label{Subsection graphs of spaces}

The fundamental groups of graphs are always free groups.  The groups we consider will usually not be free, so we need a more general notion.  A \emph{graph of spaces} $\Gamma$ consists of:
\begin{enumerate}
\item a set $V(\Gamma)$ of connected spaces, called \emph{vertex spaces};
\item a set $E(\Gamma)$ of connected spaces, called \emph{edge spaces};
\item for each edge space $e\in E(\Gamma)$ a pair of $\pi_1$-injective continuous \emph{edge maps}
$$
\partial_\pm^e:e\to\coprod_{V\in V(\Gamma)}V.
$$
When the edge in question is unambiguous, we often suppress the superscript and refer to $\partial_\pm^e$ simply as $\partial_\pm$.
\end{enumerate}

The associated topological space $|\Gamma|$ is defined as the quotient of
$$
\coprod_{V\in V(\Gamma)} V \sqcup \coprod_{e\in E(\Gamma)}(e\times [-1,+1])
$$
obtained by identifying $(x,\pm 1)$ with $\partial_\pm^e(x)$ for each edge space $e$ and every $x\in e$.  We will usually assume that $|\Gamma|$ is connected.  If $X=|\Gamma|$ we will often say that \emph{$\Gamma$ is a graph-of-spaces decomposition for $X$}, or just that \emph{$X$ is a graph of spaces}.  The \emph{underlying graph} of $\Gamma$ is the abstract graph given by replacing every vertex and edge space of $\Gamma$ by a point.  If $X$ is the topological space associated to the graph of spaces $\Gamma$, the underlying graph is denoted $\Gamma_X$. Note that there is a natural $\pi_1$-surjective map $\phi_X:X\to\Gamma_X$.

Given a graph of spaces $X$, consider a subgraph $\Gamma'\subset\Gamma_X$.  The corresponding graph of spaces $X'=\phi_X^{-1}(\Gamma')$ has a natural inclusion $X'\hookrightarrow X$ and is called a \emph{sub-graph of spaces} of $X$.  If $X$ and $Y$ are graphs of spaces, a map $f:X\to Y$ is a \emph{map of graphs of spaces} if, whenever $Y'$ is a sub-graph of spaces of $Y$, the pre-image $f^{-1}(Y')$ is a sub-graph of spaces of $X$.

\begin{defn}
A connected sub-graph of spaces $X'\subset X$ with finite underlying graph such that the inclusion map is a $\pi_1$-isomorphism is called a \emph{core} for $X$.
\end{defn}
Note that if $\pi_1(X)$ is finitely generated then a core always exists.

The fundamental group of a graph of spaces is naturally a graph of groups by the Seifert--van Kampen Theorem.  For more on graphs of spaces, graphs of groups and Bass--Serre Theory see \cite{SW} and \cite{S77}.  Our notation for graphs of spaces is different to that of \cite{SW}.

\subsection{Limit groups}

Fix a finitely generated non-abelian free group $\F$ of rank at least 2.

\begin{defn}
A group $G$ is \emph{residually free} if, for any $g\in G\smallsetminus 1$, there exists a homomorphism $f:G\to\F$ with $f(g)\neq 1$.  Likewise, $G$ is \emph{fully residually free} (or \emph{$\omega$-residually free}) if, for any finite subset $S\subset G$, there exists a homomorphism $f:G\to\F$ such that $f|_S$ is injective.  A finitely generated, fully residually free group is called a \emph{limit group}.
\end{defn}

Sela's original definition of limit groups is different, but the two definitions are equivalent (see \cite{Se1}, Definition 1.2 and Theorem 4.6).  We will need three related properties of limit groups.  A group is \emph{CSA} (standing for \emph{conjugately separated abelian}) if every maximal abelian subgroup is malnormal.

\begin{defn}
Let $\Gamma$ be a graph-of-groups decomposition for a group $G$, so $G$ acts on the Bass--Serre tree $T$.  Endow $T$ with a geodesic metric in which every edge has length 1.  Then $\Gamma$ is \emph{$k$-acylindrical} if, for every $g\in G\smallsetminus 1$, the set of fixed points of the action of $g$ on $T$ is either empty or has diameter at most $k$.  If $\Gamma$ is $k$-acylindrical for some finite $k$ then $\Gamma$ is called \emph{acylindrical}.
\end{defn}

\begin{lem}\label{Limit group properties}
Let $G$ be a limit group and consider a non-trivial decomposition $\Gamma$ of $G$ as an amalgamated free product $G=A*_C B$ with $C$ abelian.
\begin{enumerate}
\item $G$ is CSA.
\item Suppose $M\subset G$ is a non-cyclic abelian subgroup.  Then $M$ is conjugate into $A$ or $B$.
\item $\Gamma$ is 2-acylindrical.
\end{enumerate}
\end{lem}
For the proofs of these assertions, see Lemmas 1.4, 2.1 and 2.3 of \cite{Se1}.

\subsection{$\ICE$ spaces}

There is a powerful structure theory for limit groups.  Let $G'$ be a group and $Z\subset G'$ the centralizer of an element. Then the group $G=G'*_Z (Z\times \Z^n)$ is said to be obtained from $G'$ by \emph{extension of a centralizer}.  We call $Z\times \Z^n$ the \emph{extended centralizer}.

\begin{defn}
Let $\ICE$ (standing for \emph{iterated centralizer extension}) be the smallest class of groups containing all finitely generated free groups that is closed under extension of centralizers.  A group in $\ICE$ is called an $\ICE$ group.
\end{defn}
All $\ICE$ groups are limit groups (see, for example, Theorem 4.15 of \cite{Wilt2}).

\begin{rem}
Suppose $G$ is obtained from $G'\in\ICE$ by extensions of a centralizer and the centralizer $Z$ is non-cyclic.  By assertion 1 of Lemma \ref{Limit group properties}, $Z$ is abelian.  By assertion 2 of Lemma \ref{Limit group properties} and induction, $Z$ is conjugate to some previous extended centralizer $Z'\times \Z^m$. It follows that, in the construction of $\ICE$, it can always be assumed that $Z$ is cyclic.
\end{rem}
Every $G\in\ICE$ has a natural finite $K(G,1)$, constructed as follows. If $G$ is free then $K(G,1)$ can be taken to be a compact graph of suitable rank.  If $G$ is obtained from $G'$ by extension of a (cyclic) centralizer and $Y=K(G',1)$ then given an essential closed curve $\partial_+:S^1\to Y$ representing a generator of $Z$ and a coordinate circle $\partial_-:S^1\to T$, a suitable $X=K(G,1)$ can be obtained from
$$
Y\sqcup ([-1,+1]\times S^1) \sqcup T
$$
by identifying $(\pm1,\theta)\in [-1,+1]\times S^1$ with $\partial_\pm(\theta)$. We call the resulting class of spaces \emph{$\ICE$ spaces}. Note that each $\ICE$ space $X$ has an associated graph-of-spaces decomposition with vertex spaces $Y$ and $T$ and edge space a circle.

The importance of the class $\ICE$ lies in the following theorem of O.~Kharlampovich and A.~Miasnikov.

\begin{thm}[\cite{KM98b}]\label{Limit groups embed in ICE groups}
A group is a limit group if and only if it is a finitely generated subgroup of an $\ICE$ group.
\end{thm}
Another proof using different techniques is given in \cite{CG04}.

Our strategy to prove Theorems \ref{Theorem A} and \ref{Theorem B} is to prove them by induction for all $\ICE$ groups, since subgroup separability and local retractions both pass to subgroups.  As mentioned in subsection \ref{Motivating example subsection}, all previous combination theorems for subgroup separability rely on some of the vertex groups being free.  This is not true of $\ICE$ groups.  Therefore we will need to prove a strong result, analogous to Corollary \ref{Motivating example}, for all groups in $\ICE$.

\section{Elevations and pre-covers}

\subsection{Elevations}

Let $X$ be a graph of spaces and $X'\to X$ a covering map. Then $X'$ inherits a graph-of-spaces decomposition. We want to understand the edge maps of $X'$.  To do this we study \emph{elevations}, natural generalizations of lifts that were introduced by Wise in \cite{Wise1}.

\begin{defn}\label{Elevation definition}
Consider a continuous map of connected based spaces $f:(A,a)\to (B,b)$ and a covering $B'\to B$.  An \emph{elevation} of $f$ to $B'$ consists of a connected covering $p:(A',a')\to (A,a)$ and a lift $f':A'\to B'$ of  $f\circ p$ so that for every intermediate covering
$$
(A',a')\to(\bar{A},\bar{a})\stackrel{q}{\to}(A,a)
$$
there is no lift $\bar{f}$ of $f\circ q$ to $B'$ with $\bar{f}(\bar{a})=f'(a')$.

Two elevations $f'_1:A'_1\to B'$ and $f'_2:A'_2\to B'$ are \emph{isomorphic} if there exists a homeomorphism $\iota:A'_1\to A'_2$, covering the identity map on $A$, such that
$$
f'_1=f'_2\circ\iota.
$$
\end{defn}

\begin{figure}[ht]
\begin{center}
\includegraphics[width=0.7\textwidth]{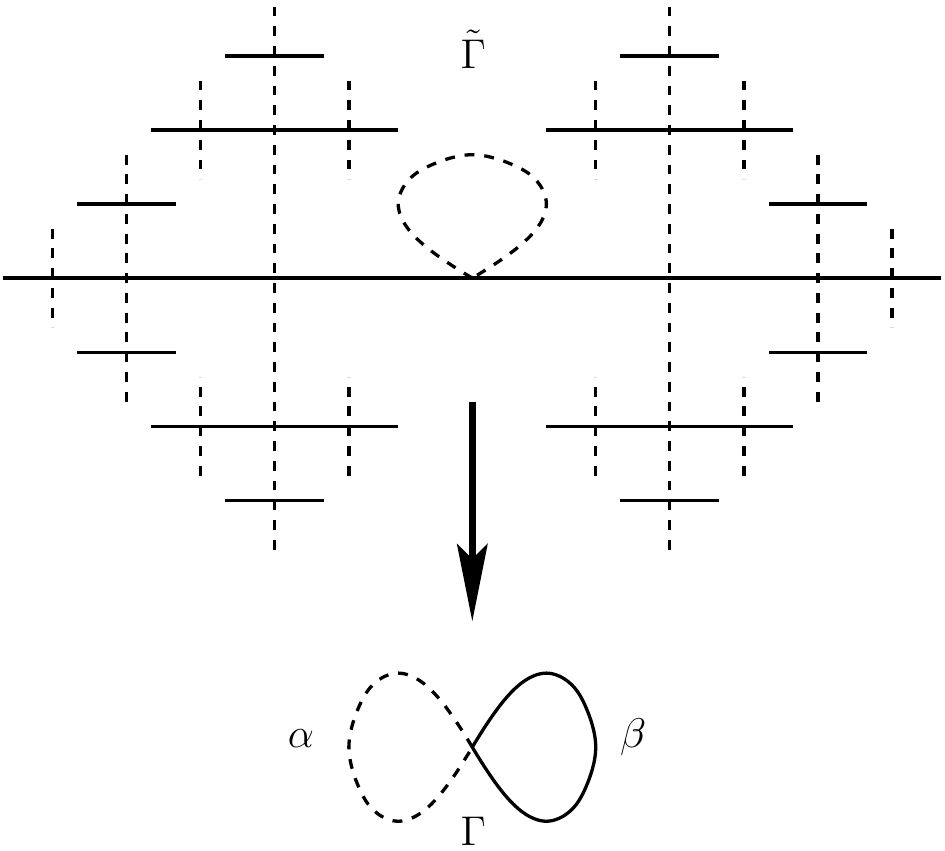}
\caption{The diagram shows a graph $\Gamma$ that is the wedge of two loops---$\alpha$ (which is dashed) and $\beta$---and a covering space $\tilde{\Gamma}$.  The dashed lines are the images of the elevations of $\alpha$ to $\tilde{\Gamma}$ and the solid lines are the images of the elevations of $\beta$ to $\tilde{\Gamma}$.  One elevation of $\alpha$ is of finite degree.  All the other elevations are of infinite degree.}
\label{Elevation to cover figure}
\end{center}
\end{figure}

A simple example of some elevations is illustrated in Figure \ref{Elevation to cover figure}.  In practice, we will often abuse notation and refer to just the lift $f'$ as an elevation of $f$.

\begin{defn}[Degree of an elevation]\label{Definition of degree}
The \emph{degree} of the elevation $f':A'\to  B'$, denoted $\deg(A')$, is the \emph{conjugacy class} of the subgroup $\pi_1(A')\subset\pi_1(A)$.  If $\pi_1(A')$ is of finite index in $\pi_1(A)$ then $f'$ is called \emph{of finite degree}; otherwise, $f'$ is \emph{of infinite degree}.
\end{defn}

\begin{rem}\label{Elevations of homotopic maps}
Consider a map $f_0:(A,a)\to (B,b_0)$ and an elevation $f'_0:(A',a')\to (B',b'_0)$.  Let $h:A\times[0,1]\to B$ be a free (unbased) homotopy taking $f_0$ to $f_1:(A,a)\to (B,b_1)$. Then $h$ lifts to a homotopy $h':A'\times [0,1]\to B'$ taking $f'_0$ to $f'_1:(A',a')\to (B',b'_1)$.  Now $f'_1$ is an elevation of $f_1$.
\end{rem}
In the light of Remark \ref{Elevations of homotopic maps}, it makes sense to consider elevations of a map defined up to free homotopy.

\begin{rem}
If $X$ has a graph-of-spaces decomposition $\Gamma$ and $X'\to X$ is a covering space then $X'$ inherits a graph-of-spaces decomposition $\Gamma_{X'}$, with vertex spaces the connected components of the pre-images of the vertex spaces of $X$ and edge spaces and maps given by all the (isomorphism classes of) elevations of the edge maps to the vertex spaces of $X'$.
\end{rem}

\begin{rem}
If $f'_1:A'\to B'$ and $f'_2:A'\to B'$ are both elevations of $f:A\to B$ and $f'_1(a')=f'_2(a')$ for some $a'\in A'$ then $f'_1$ and $f'_2$ are isomorphic.
\end{rem}

Elevations can also be understood algebraically, in terms of fundamental groups.

\begin{lem}\label{Covering space theory}
Let $f:(A,a)\to(B,b)$ and $B'\to B$ be as in Definition \ref{Elevation definition}.  Fix a lift $b'\in B'$ of $b=f(a)$, so $\pi_1(B',b')$ is identified with a subgroup of $\pi_1(B,b)$.  Consider the covering $p:(A',a')\to(A,a)$ such that
$$
\pi_1(A',a')=f_*^{-1}(\pi_1(B',b')).
$$
The composition $f\circ p$ admits a lift $f':A'\to B'$ to $B'$, which is an elevation of $f$.  Every elevation of $f$ to $B'$ arises in this way.
\end{lem}
The proof is standard covering-space theory. See, for example, Proposition 1.33 of \cite{Ha02}.

Fix $g\in\pi_1(B,b)$. The terminus of the lift at $b'$ of $g$ to $B'$ gives a well-defined new basepoint $\beta\in B'$, and so defines an elevation $f':(A',\alpha)\to (B',\beta)$ of $f$ with
$$
\pi_1(A',\alpha)= f_*^{-1}\pi_1(B',\beta).
$$
Note that $\pi_1(B',\beta)=\pi_1(B',b')^g$.  When do two elements of $\pi_1(B,b)$ yield isomorphic elevations?

\begin{lem}[Elevations correspond to double cosets]\label{Double cosets}
Consider the situation of Lemma \ref{Covering space theory}.  Let $g_1,g_2\in\pi_1(B,b)$ and let $f'_i:(A'_i,\alpha_i)\to (B',\beta_i)$ be the elevation of $f$ such that
$$
\pi_1(A'_i,\alpha_i)=f_*^{-1}(\pi_1(B',b')^{g_i}).
$$
Then $f'_1$ and $f'_2$ are isomorphic elevations if and only if
$$
\pi_1(B',b')g_1(f_*\pi_1(A,a))=\pi_1(B',b')g_2(f_*\pi_1(A,a))
$$
in $\pi_1(B',b')\backslash \pi_1(B,b)/(f_*\pi_1(A,a))$.
\end{lem}
\begin{pf}
Suppose $f'_1$ and $f'_2$ are isomorphic, so there exists a covering transformation $\iota:A'_1\to A'_2$ such that $f'_1=f'_2\circ\iota$.  Let $\gamma':I\to A'_2$ be a path from $\alpha_2$ to $\iota(\alpha_1)$.  Fix representative loops for $g_1$ and $g_2$, and let $g'_1$ and $g'_2$ be their respective lifts at $b'$ to $B'$. The concatenation $g'_2\cdot f'_2\circ\gamma'\cdot (g'_1)^{-1}$ makes sense and determines an element of $\pi_1(B',b')$.  Since the image $\gamma$ of $\gamma'$ in $A$ determines an element of $\pi_1(A,a)$, we have $g_2f_*(\gamma)g_1^{-1}\in\pi_1(B',b')$.  It follows that
$$
\pi_1(B',b')g_1(f_*\pi_1(A,a))=\pi_1(B',b')g_2(f_*\pi_1(A,a))
$$
as required.

Now suppose that
$$
\pi_1(B',b')g_1(f_*\pi_1(A,a))=\pi_1(B',b')g_2(f_*\pi_1(A,a)).
$$
Then there exists $\gamma\in\pi_1(A,a)$ such that $g_2f_*(\gamma)g_1^{-1}\in\pi_1(B',b')$ and so
$$
\pi_1(B',b')^{g_2f_*(\gamma)}=\pi_1(B',b')^{g_1}.
$$
Since $\pi_1(A'_i,\alpha_i)=f_*^{-1}(\pi_1(B',b')^{g_i})$,
$$
\pi_1(A'_2,\alpha_2)^\gamma=\pi_1(A'_1,\alpha_1).
$$
Let $\gamma'$ be the lift at $\alpha_2$ of $\gamma$ to $A'_2$.  By standard covering space theory, there exists a covering transformation $\iota:A'_1\to A'_2$ so that $\iota(\alpha_1)$ is the terminus of $\gamma'$. Since $f'_2\circ\iota$ and $f'_1$ are elevations of $f$ that agree at $\alpha_1$ it follows that $f'_2\circ\iota=f'_1$.
\end{pf}

Henceforth, we will often suppress mention of basepoints.   To summarize the above discussion, the (isomorphism classes of) elevations of a map $f:A\to B$ to a cover $B'$ are identified with the set of double cosets
$$
\pi_1(B')\backslash\pi_1(B)/(f_*\pi_1(A)).
$$

Consider a map $f:A\to B$, a connected covering $B'\to B$ and an elevation $f':A'\to B'$ of $f$.  Let
$$
B'\to\bar{B}\to B
$$
be an intermediate covering.  Then there exists a unique elevation $\bar{f}:\bar{A}\to\bar{B}$ of $f$ such that $A'\to A$ factors through $\bar{B}\to B$ and the diagram

\[ \begin{CD}
{A'} @>f'>>  {B'} \\
@VVV    @VVV \\
 {\bar{A}} @>\bar{f}>>  {\bar{B}}\\
\end{CD}\]
commutes, determined by the requirement that $\pi_1(\bar{A})=f_*^{-1}(\pi_1(\bar{B}))$.

\begin{defn}
In the situation of the previous paragraph, we say $f'$ \emph{descends to $\bar{f}$}.
\end{defn}

\begin{rem}
If $q:B'\to \bar{B}$ is injective on $f'(A')$ then $f'$ descends to the composition $\bar{f}=q\circ f'$; in particular, $\deg{\bar{f}}=\deg{f'}$. Furthermore, if $f'_1:A'_1\to B'$ and $f'_2:A'_2\to B'$ are non-isomorphic elevations of $f$ descending to $\bar{f}_1$ and $\bar{f}_2$ respectively, and $q$ is injective on $f'_1(A'_1)\cup f'_2(A_2)$, then $\bar{f}_1$ and $\bar{f_2}$ are non-isomorphic.
\end{rem}

As one might expect, elevations to normal covers behave well.

\begin{lem}\label{Elevations to normal covers}
Consider a map $f:(A,a)\to (B,b)$ and a normal covering
$$
(\tilde{B},\tilde{b})\to (B,b).
$$
For any elevation
$$
\tilde{f}:(\tilde{A},\tilde{a})\to(\tilde{B},\tilde{b})
$$
of $f$ to $\tilde{B}$:
\begin{enumerate}
\item $\tilde{A}$ is a normal cover of $A$;
\item $f_*$ descends to a monomorphism
$$
\tilde{f}_*:\pi_1(A)/\pi_1(\tilde{A})\to\pi_1(B)/\pi_1(\tilde{B});
$$
\item and $\tilde{f}$ is $\pi_1(A)$-equivariant---that is, for $g\in\pi_1(A)$,
$$
\tilde{f}\circ\tau_g=\tau_{\tilde{f}_*(g)}\circ\tilde{f}
$$
where $\tau_g$ denotes the covering transformation by which $g$ acts on $\tilde{A}$ and, likewise, $\tilde{f}_*(g)$ acts on $\tilde{B}$ by $\tau_{\tilde{f}_*(g)}$.
\end{enumerate}
\end{lem}
\begin{pf}
Since $\pi_1(\tilde{A})=f_*^{-1}\pi_1(\tilde{B})$ is a normal subgroup of $\pi_1(A)$, $\tilde{A}$ is a normal cover.  Since $\pi_1(\tilde{A})=f_*^{-1}\pi_1(\tilde{B})$, $f_*$ descends to an injective map $\pi_1(A)/\pi_1(\tilde{A})\to\pi_1(B)/\pi_1(\tilde{B})$. Both $\tilde{f}\circ\tau_g$ and $\tau_{\tilde{f}_*(g)}\circ\tilde{f}$ are elevations of $f$ that map $\tilde{a}$ to the terminus of the lift at $\tilde{b}$ of $f_*(g)$ to $\tilde{B}$. So they are equal.
\end{pf}

\subsection{Axes of translation}

As remarked in subsection \ref{Subsection graphs of spaces}, a graph of spaces $X$ is equipped with a map to the underlying graph $\phi_X:X\to\Gamma_X$. One makes $\Gamma_X$ into a graph of groups by labelling each vertex and edge with the fundamental group of its pre-image.  If $\tilde{X}\to X$ is the universal covering then the underlying graph $\Gamma_{\tilde{X}}$ is naturally identified with the Bass--Serre tree of $\Gamma_X$ and the map $\phi_{\tilde{X}}:\tilde{X}\to\Gamma_{\tilde{X}}$ is $\pi_1(X)$-equivariant.  (There is a subtlety here to do with basepoints.  The action on $\tilde{X}$ and the action on $\Gamma_{\tilde{X}}$ both depend on choices of basepoints.  The map $\phi_{\tilde{X}}$ is equivariant if it maps the basepoint of $\tilde{X}$ to the basepoint of $\Gamma_{\tilde{X}}$.)  More generally, if $X^H\to X$ corresponds to the subgroup $H\subset \pi_1(X)$ then the underlying graph $\Gamma_{X^H}$ is naturally identified with the quotient graph $H\backslash\Gamma_{\tilde{X}}$.

It is a fundamental fact about the action of a group $G$ on a tree that every $g\in G$ either fixes a point or acts by translating an embedded line, denoted $\Axis(g)$ (see, for instance, \cite{S77}). Elements of the first sort are called elliptic and elements of the second sort are called hyperbolic. Accordingly, a loop $\delta:S^1\to X$ is called \emph{elliptic} if it is freely homotopic into a vertex space; otherwise, it is called \emph{hyperbolic}.

\begin{lem}\label{Axis in elevation}
Let $(X,x_0)$ be a (based) graph of spaces and let $(\tilde{X},\tilde{x}_0)$ be the (based) universal cover.  Let $\delta:S^1\to X$ be a (based) loop and let $\tilde{\delta}:\tilde{S}^1\to\tilde{X}$ be the elevation to $\tilde{X}$ corresponding to the coset $g\langle\delta\rangle$.  Then $\tilde{\delta}(\tilde{S}^1)$ is a $\tau_{\delta^{(g^{-1})}}$-invariant subspace (using the notation of Lemma \ref{Elevations to normal covers}).  Therefore:
\begin{enumerate}
\item if $\delta$ is elliptic then $\phi_{\tilde{X}}\circ\tilde{\delta}(\tilde{S}^1)$ contains a $\delta^{(g^{-1})}$-invariant vertex;
\item if $\delta$ is hyperbolic then $\phi_{\tilde{X}}\circ\tilde{\delta}(\tilde{S}^1)$ contains $\Axis(\delta^{(g^{-1})})$.
\end{enumerate}
\end{lem}
\begin{pf}
Let $\tau_g$ be the covering transformation of $\tilde{X}$ corresponding to $g$.  Since $\tilde{\delta}$ is the elevation of $\delta$ at the terminus of the lift of $g$ at $\tilde{x}_0$ it follows that $\tau_g^{-1}\circ\tilde{\delta}$ is the elevation of $\delta$ at $\tilde{x}_0$.  By Lemma \ref{Elevations to normal covers}, therefore, $\tau_g^{-1}\circ\tilde{\delta}$ is a $\delta$-equivariant map.  Hence $\tilde{\delta}$ is a $\delta^{(g^{-1})}$-equivariant map.

As $\phi_{\tilde{X}}:\tilde{X}\to\Gamma_{\tilde{X}}$ is also equivariant we have that $\phi_{\tilde{X}}\circ\tilde{\delta}(\tilde{S}^1)$ is $\delta^{(g^{-1})}$-invariant.  It is a fact about group actions on trees that any subtree invariant under the action of a hyperbolic element contains the axis and any subtree invariant under the action of an elliptic element contains a fixed point.
\end{pf}

\subsection{Acylindrical graphs of groups and proper loops}

The fact that an ICE group decomposes as an acylindrical graph of groups has profound implications for the structure of limit groups.

\begin{defn}\label{Proper loops defn}
Let $X$ be a graph of spaces and consider an arbitrary covering $X'\to X$ with $\pi_1(X')$ finitely generated; let $\phi':X'\to\Gamma_{X'}$ be the map to the underlying graph.  A loop
$$
\gamma:C\cong S^1\to X
$$
is called \emph{proper} if, for any such covering $X'$ and any elevation $\gamma'$ of $\gamma$ to $X'$, the composition map
$$
\lambda=\phi'\circ\gamma':C'\to \Gamma_{X'}
$$
is a proper map.
\end{defn}

In the case when $X$ really is a graph, all loops are proper.

\begin{lem}\label{Loops in graphs are proper}
Let $X$ be a compact graph.  Then any loop $\gamma:C\to X$ is proper.
\end{lem}
\begin{pf}
Let $X'\to X$ be a covering with $H=\pi_(X')$ finitely generated and let $\gamma':C'\to X'$ be an elevation.  If $C'$ is compact then there is nothing to prove, so assume that $C'\cong\R$.  Let $\tilde{X}$ be the universal cover of $X$ and let $\tilde{X}_H$ be the minimal $H$-invariant subtree.  If $\gamma'$ is not proper then $\Axis(\gamma)$ intersects $\tilde{X}_H$ in infinitely many vertices and therefore an infinite ray.  Consider a vertex $v$ of $\tilde{X}_H$.  There are only finitely many $H$-orbits of vertices in $\tilde{X}_H$ so, for some integer $n$ and some element $h\in H$, $\gamma^nh^{-1}$ fixes $v$.  As the action on $\tilde{X}$ is free it follows that $\gamma^n\in H$, so $C'$ is compact, contradicting the assumption that $C'\cong\R$.
\end{pf}

This lemma is implicit in the proof of Corollary \ref{Motivating example}.  We will need a similar result for $\ICE$ spaces, which we will prove by induction.

\begin{lem}
Let $\Gamma$ be a $k$-acylindrical graph of groups and consider the action of $G=\pi_1(\Gamma)$ on the Bass--Serre tree $\tilde{\Gamma}$.  Suppose that $\gamma\in G$ is a hyperbolic element with translation length $|\gamma|$.  Let $g$ be another element of $G$.  If the intersection $\Axis(\gamma)\cap g\Axis(\gamma)$ is of length greater than $k+|\gamma|$ then $g$ commutes with $\gamma$.  In particular, $g\Axis(\gamma)=\Axis(\gamma)$.
\end{lem}
\begin{pf}
Since the axis of $\gamma$ can be defined as the set of points that $\gamma$ moves a minimal distance, $g\Axis(\gamma)=\Axis(\gamma^{g^{-1}})$ and $|\gamma|=|\gamma^{g^{-1}}|$.  If the length of $\Axis(\gamma)\cap g\Axis(\gamma)$ is greater than $k+|\gamma|$ then there exists an arc $\alpha$ of length greater than $k$ such that $\gamma^{g^{-1}}\gamma^{-1}$ fixes $\alpha$ pointwise.  Since the action is $k$-acylindrical, $\gamma^{g^{-1}}\gamma^{-1}=[g,\gamma]=1$.
\end{pf}

We use this fundamental fact in the following technical lemma. For points $a,b\in \tilde{\Gamma}$, denote by $[a,b]$ the unique geodesic arc from $a$ to $b$.

\begin{lem}
Let $X$ be a graph of spaces for which the corresponding graph of groups is $k$-acylindrical.  Let $\gamma\in\pi_1(X)$ be a proper loop.  Then for any covering $X'\to X$ with finitely generated fundamental group and for any non-trivial $g\in \pi_1(X)$, there are at most finitely many cosets $\gamma^n\pi_1(X')$ in $\pi_1(X)/\pi_1(X')$ such that
\[
\gamma^{-n} g\gamma^n\in\pi_1(X').
\]
\end{lem}
\begin{pf}
Let $G=\pi_1(X)$ and let $H=\pi_1(X')$.  Without loss of generality, $g\in H$.  If $\langle\gamma\rangle\cap\pi_1(X')\neq 1$ then there is nothing to prove, so assume that $\langle\gamma\rangle\cap\pi_1(X') = 1$.  Let $\tilde{\Gamma}$ be the Bass--Serre tree of the induced splitting of $G$ equipped with the standard combinatorial metric $d$ and let $\tilde{\Gamma}_H$ be the minimal $H$-invariant subtree.  Since $\gamma$ is proper, $\tilde{\Gamma}_H\cap\Axis(\gamma)$ is a compact arc $\alpha$.  Any element $h\in H$ that commutes with $\gamma$ must fix $\alpha$, and so the whole of $\Axis(\gamma)$, pointwise.  So, as the splitting is acylindrical, any such $h$ is trivial.

Fix a basepoint $t_0\in\alpha$.  Let $|\gamma|$ be the translation length of $\gamma$ and let $\beta$ be the subarc of $\Axis(\gamma)$ consisting of those points that are a distance at most $k+|\gamma|$ from $\tilde{\Gamma}_H$.  Let $n$ be large enough that $\gamma^{\pm n}(\beta)\cap \tilde{\Gamma}_H=\varnothing$.

Let $s$ be the nearest point on $\Axis(\gamma)$ to $g\gamma^nt_0$.  Note that $g\gamma^nt_0\in g\Axis(\gamma)$.  Since $g$ is non-trivial it does not commute with $\gamma$, so the diameter of $\Axis(\gamma)\cap g\Axis(\gamma)$ is at most $k+|\gamma|$ and hence $s\in\beta$. Therefore $\gamma^{-n}s\notin \tilde{\Gamma}_H$. But the geodesic arc $[g\gamma^n t_0,\gamma^{-n}g\gamma^n t_0]$ intersects $\Axis(\gamma)$ precisely in the arc $[s,\gamma^{-n}s]$.  It follows that $\gamma^{-n}g\gamma^nt_0\notin\tilde{\Gamma}_H$ so $\gamma^{-n}g\gamma^n\notin H$.
\end{pf}

Using this observation, we can prove the extremely useful fact that all hyperbolic loops in ICE spaces are proper.

\begin{lem}[Hyperbolic loops in $\ICE$ spaces are proper]\label{Proper loops}
Let $X$ be an $\mathrm{ICE}$ space.  Every hyperbolic loop $\delta:C\to X$ is proper.
\end{lem}
\begin{pf}
As usual $X$ is constructed by gluing a torus $T$ to an $\mathrm{ICE}$ space $Y$ of lower level.  By induction on level it can be assumed that the lemma holds for $Y$.  Without loss of generality take $\delta$ to be based.  The graph-of-spaces decomposition $\Gamma$ for $X$ induces a graph-of-groups decomposition of $\pi_1(X)$ which, abusing notation, we also denote by $\Gamma$.  Decompose $\delta$ in $\Gamma$ as an irreducible word
$$
\delta=a_0b_1a_1\ldots b_na_n
$$
where $a_i\in\pi_1(T)$ and $b_i\in\pi_1(Y)$.  Let $X'\to X$ be a covering of $X$ and let $\delta':C'\to X'$ be an elevation of $\delta$, as in Definition \ref{Proper loops defn}.

If $C'$ is compact then there is nothing to prove, so assume $C'\cong\R$.  Let $\pi_1(e)=\langle z\rangle$.  If $\lambda=\phi'\circ\delta'$ is not proper then for some vertex space $Y'$ of $X'$ and two elevations $\partial^+_{1,2}:e'_{1,2}\to Y'$ the line $\lambda$
crosses $Y'$ from $\partial^+_1(e'_1)$ to $\partial^+_2(e'_2)$ in infinitely many places.  Therefore, without  loss of generality, for infinitely many indices $p$,
non-zero indices $q_p$ and some $b_k$,
$$
z^p [b_k,z^{q_p}]z^{-p} \in \pi_1(Y').
$$
Since the elevations $\partial^+_{1,2}$ are proper by induction, it follows that $[b_k,z]=1$, contradicting the irreducibility of the decomposition of $\delta$.
\end{pf}

\subsection{Pre-covers}

Pre-covers fill the role played by immersions of graphs in Stallings' proof of Hall's Theorem.  In \cite{Gi1}, Gitik uses a notion of pre-cover; our pre-covers are analogous to hers.

\begin{figure}[htp]
\begin{center}
\psfrag{Gamma'}{$\bar{\Gamma}$}\psfrag{Gammatilde}{$\tilde{\Gamma}$}
\includegraphics[width=0.7\textwidth]{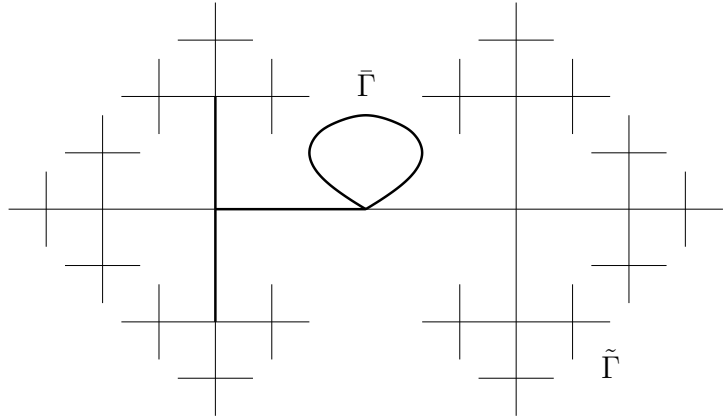}
\caption{The diagram shows the covering space $\tilde{\Gamma}$ from Figure \ref{Elevation to cover figure}.  A finite subgraph $\bar{\Gamma}$ of $\tilde{\Gamma}$ is drawn with a thick line. The restriction of the covering map $\tilde{\Gamma}\to\Gamma$ is a pre-covering.  The canonical completion of $\bar{\Gamma}$ is $\tilde{\Gamma}$, which can be constructed by attaching the forest drawn with thin lines.}\label{Pre-cover figure}
\end{center}
\end{figure}

\begin{defn}
Let $X$ and $X'$ be graphs of spaces ($X'$ is not assumed to be connected). A \emph{pre-covering} is a locally injective map $X'\to X$ that maps vertex spaces and edge spaces of $X'$ to vertex spaces and edge spaces of $X$ respectively, and restricts to a covering on each vertex space and each edge space. Furthermore, for each edge space $e'$ of $X'$ mapping to an edge space $e$ of $X$, the diagram of edge maps
\[ \begin{CD}
{e'} @>{\bar{\partial}_\pm}>>  {V'_\pm} \\
@VVV    @VVV \\
 {e} @>{\partial_\pm}>>  {V_\pm}\\
\end{CD}\]
is required to commute.  The domain $X'$ is called a \emph{pre-cover}.  An example of a pre-cover, in the context of graphs, is illustrated in Figure \ref{Pre-cover figure}.

The pre-covering $X'\to X$ is \emph{finite-sheeted} if the pre-image of every point of $X$ is finite.
\end{defn}

\begin{rem}
All the edge maps of $X'$ are elevations of edge maps of $X$ to the vertex spaces of $X'$.  An elevation of an edge map of $X$ to a vertex space of $X'$ that is not (isomorphic to) an edge map of $X'$ is called \emph{hanging}.  If none of the elevations are hanging then $X'$ is in fact a cover.
\end{rem}

\begin{prop}[Canonical completion of pre-covers]\label{Canonical completion of pre-covers}
Let $X$ be a graph of spaces and $X'\to X$ a pre-covering with $X'$ connected. Then $\pi_1(X')$ injects into $\pi_1(X)$ and, furthermore, there exists a unique embedding $X'\hookrightarrow\tilde{X}$ into a connected covering $\tilde{X}\to X$ such that $\pi_1(\tilde{X})=\pi_1(X')$ and $\tilde{X}\to X$ extends $X'\to X$. This is called the \emph{canonical completion} of $X'\to X$.
\end{prop}
For the proof, see Proposition 2.9 of \cite{Wilt1}.  In this situation, $X'$ is a core for $\tilde{X}$.  It should be stressed that this construction does not preserve finiteness of sheets.  Indeed, the completion $\tilde{X}$ is finite-sheeted if and only if $\pi_1(X')$ is of finite index in $\pi_1(X)$.

\begin{figure}[htp]
\begin{center}
\psfrag{Gamma'}{$\bar{\Gamma}$}\psfrag{Gammatilde}{$\tilde{\Gamma}$}\psfrag{alpha1}{$\tilde{\alpha}_1$}\psfrag{alpha2}{$\tilde{\alpha}_2$}
\includegraphics[width=0.7\textwidth]{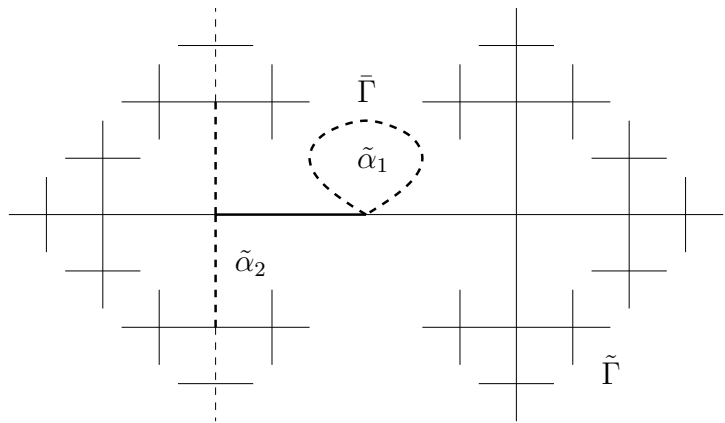}
\caption{Two elevations, $\tilde{\alpha}_1$ and $\tilde{\alpha}_2$, of $\alpha$ to $\tilde{\Gamma}$ are marked with dashed lines.   Their restrictions to $\bar{\Gamma}$, $\bar{\alpha}_1$ and $\bar{\alpha}_2$, are examples of elevations of $\alpha$ to $\bar{\Gamma}$.  The elevation $\bar{\alpha}_1$ is full whereas $\bar{\alpha}_2$ is not.  Both elevations are disparate.  Likewise, the set of elevations $\{\bar{\alpha}_1,\bar{\alpha}_2\}$ is also disparate.}  \label{Elevations to pre-covers figure}
\end{center}
\end{figure}

\begin{defn}[Elevations to pre-covers]\label{Elevations to pre-covers}
Let $f:A\to B$ be a map of graphs of spaces and $B'\to B$ a pre-covering.  Let $\tilde{B}\to B$ be the canonical completion of $B'\to B$. Consider an elevation $\tilde{f}:\tilde{A}\to\tilde{B}$ of $f$. Let $A'=\tilde{f}^{-1}(B')\subset \tilde{A}$. Note that $A'$ is a pre-cover of $A$ since $f$ is a map of graphs of spaces. If $A'$ is non-empty then
$$
f'=\tilde{f}|_{A'}:A'\to B'
$$
is called an \emph{elevation of $f$ to $B'$}.

When $A'\to A$ is a genuine covering of $A$ the elevation $f'$ is called \emph{full}.  In this case, it makes sense to define the \emph{degree} of $f'$ to be the conjugacy class of $\pi_1(A')$ in $\pi_1(A)$.
\end{defn}

Note that this definition of an elevation to a pre-cover is slightly less general than the definition given in \cite{Wilt1}. However, it is less complicated and will suffice for the purposes of this paper.  Two examples of elevations to a pre-cover are illustrated in Figure \ref{Elevations to pre-covers figure}.

\subsection{Stallings' Principle}\label{Subsection Stallings'
Principle}

With the language of pre-covers and elevations we can generalize Theorem \ref{Immersions of graphs} to extend finite-sheeted pre-covers of graphs of spaces to finite-sheeted covers. We call this generalization \emph{Stallings' Principle}.  Let $X$ be a graph of spaces and $p:\bar{X}\to X$ a pre-covering.   Consider an edge map $\partial:e\to V$ of $X$ and an elevation $\bar{\partial}:\bar{e}\to\bar{V}$ of $\partial$.  Since the image of $\partial$ is contained in $V$ and $p$ restricts to a covering map on $\bar{V}$ it follows that $\bar{e}\to e$ is a covering so $\bar{\partial}$ is full.  Indeed, $\bar{\partial}$ is a genuine elevation to $\bar{V}$.  It therefore makes sense to consider the degree of $\bar{\partial}$ (see Definition \ref{Elevations to pre-covers}).

\begin{thm}[Stallings' Principle]\label{Stallings}
Let $\bar{X}\to X$ be a pre-covering with the following property: for each edge space $e$ of $X$ with edge maps $\partial_\pm:e\to V_\pm$, for each conjugacy class $\mathcal{D}$ of subgroups of $\pi_1(e)$, there exists a bijection between the set of elevations of $\partial_+$ to $\bar{X}$ of degree $\mathcal{D}$ and the set of elevations of $\partial_-$ to $\bar{X}$ of degree $\mathcal{D}$.  Then $\bar{X}\to X$ can be extended to a covering $\hat{X}\to X$, with the same set of vertex spaces as $\bar{X}$.
\end{thm}

For the proof, see Proposition 3.1 of \cite{Wilt1}.  The crucial feature of this construction (as opposed to, say, Proposition \ref{Canonical completion of pre-covers}) is that it extends finite-sheeted pre-covers to finite-sheeted covers: since $\hat{X}$ has the same vertex spaces as $\bar{X}$, if $\bar{X}\to X$ is finite-sheeted then so is $\hat{X}\to X$.  From Theorem \ref{Stallings} one deduces Theorem \ref{Immersions of graphs} as follows.

\begin{pfof}{Theorem \ref{Immersions of graphs}}
The immersion of graphs $X'\to X$ is a pre-covering. Let $N$ be the maximum number of pre-images in $X'$ of a vertex of $X$. Create a new (disconnected) pre-cover $\bar{X}$ by adding disconnected vertices to $X'$ so that every vertex of $X$ has $N$ pre-images. Now every edge map of $\Gamma$ has $N$ elevations to the pre-cover, all of the same (trivial) degree; so, by Theorem \ref{Stallings}, $\bar{X}\to X$ can be extended to a covering $\hat{X}\to X$.
\end{pfof}

The proof of Corollary \ref{Hall's Theorem for free groups} makes use of the fact that a covering of a graph that corresponds to a finitely generated subgroup has a finite core.  In contrast, in the case of a graph of spaces $X$, a covering $X^H$ corresponding to finitely generated subgroup $H\subset\pi_1(X)$ has a core $X'$ with finite underlying graph, but in general $X'\to X$ is not finite-sheeted.  Most of the rest of this paper is devoted to addressing this difficulty.

\subsection{Making elevations full}

Let $\delta:S^1\to X$ be a loop in $X$ and $\bar{X}\to X$ a pre-covering. Let $\bar{\delta}:\bar{S}^1\to\bar{X}$ be an elevation of $\delta$ to $\bar{X}$.  The key technical construction enables us, under the hypothesis of \emph{disparity}, to complete $\bar{\delta}$ to a full elevation---that is, an elevation whose domain is homeomorphic to $S^1$ or $\R$.  Indeed, we can do this for sets of elevations.

First we fix some notation.  Let $\tilde{X}$ be the canonical completion of $\bar{X}$ so $\bar{\delta}$ is a restriction of an elevation $\tilde{\delta}:\tilde{S}^1\to\tilde{X}$.  Given $x\in\partial\bar{S}^1$, the closure of the adjacent component of $\tilde{S}^1\smallsetminus\bar{S}^1$ can be identified with some subinterval of $[0,\infty)$, with $x$ identified with $0$.  For sufficiently small $\epsilon>0$, $\tilde{\delta}((0,\epsilon))$ is contained in a single edge space of $\tilde{X}$, corresponding to a hanging elevation $\bar{\partial}_x:e_x\to\bar{X}$ of some edge map $\partial$.

\begin{defn}[Disparity]
Let $\{\delta_i:C_i\to X\}$ be a set of loops.  Let $\bar{X}\to X$ be a pre-cover, as above.  A collection of elevations
$$
\{\bar{\delta}_j:\bar{C}_j\to\bar{X}\},
$$
where each $\bar{\delta_j}$ is an elevation of some $\delta_i$, is \emph{disparate} if:
\begin{enumerate}
\item $\coprod_j\bar{C}_j$ is compact;
\item the map from $\coprod_j\partial\bar{C}_j$ to the set of hanging elevations of edge maps given by
$$
x\mapsto\bar{\partial}_x
$$
is injective; and,
\item for each $x\in\coprod_j\partial\bar{C}_j$, $e_x$ is simply connected.
\end{enumerate}
\end{defn}

If a path $\delta$ cannot be homotoped (relative to its endpoints) off any of the edge or vertex spaces of $X$ that it intersects, then $\delta$ is called \emph{reduced}.  Clearly, any path is homotopic to a reduced path.  Recall that an elevation  $\bar{\delta}_j$ to a pre-cover is full if its domain $\bar{C}_j$ is a covering space of $C_i$.  The next lemma asserts that, under the hypothesis of disparity, we can extend a pre-cover to make a collection of elevations of reduced loops full.

\begin{figure}[htp]
\begin{center}
\psfrag{Gamma'}{$\bar{\Gamma}$}\psfrag{Gammahat}{$\hat{\Gamma}$}\psfrag{alpha1}{$\bar{\alpha}_1$}\psfrag{alpha2}{$\bar{\alpha}_2$}\psfrag{alphahat1}{$\hat{\alpha}_1$}\psfrag{alphahat2}{$\hat{\alpha}_2$}
\includegraphics[width=0.7\textwidth]{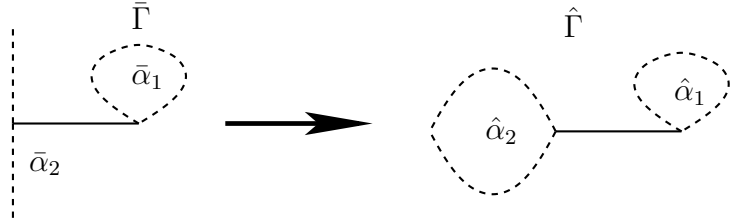}
\caption{Gluing gives an intermediate pre-cover $\hat{\Gamma}$ in which $\bar{\alpha}_2$ descends to a full elevation $\hat{\alpha}_2$.}\label{Making elevations full figure}
\end{center}
\end{figure}

\begin{lem}[Making elevations full]\label{Making elevations full}
Let $\bar{X}\to X$ be a pre-covering, let $\{\delta_i:C_i\to X\}$ be a set of reduced loops and let $\{\bar{\delta}_j:\bar{C}_j\to\bar{X}\}$ be a disparate set of elevations to $\bar{X}$.  Suppose that, for each $j$, $\hat{C}_j\to C_i$ is a finite-sheeted covering extending $\bar{C}_j\to C_i$. Then there exists a pre-covering $\hat{X}\to X$ extending $\bar{X}\to X$ so that each $\bar{\delta}_j$ extends to a full elevation
$$
\hat{\delta}_j:\hat{C}_j\to\hat{X}.
$$
If $\bar{X}$ has finite underlying graph, each $\hat{C}_j\to C_i$ is finite-sheeted and there are finitely many $\bar{\delta}_j$ then the resulting $\hat{X}$ has finite underlying graph.  Furthermore, $\pi_1(\hat{X})=\pi_1(\bar{X})*F$ where $F$ is some finitely generated free group.
\end{lem}
The proof of Lemma \ref{Making elevations full} is the same as that of Proposition 2.13 in \cite{Wilt1}.  However, because of its fundamental importance to our argument, we repeat the proof here.  Figure \ref{Making elevations full figure} illustrates the idea of Lemma \ref{Making elevations full} in the situation of Figure \ref{Elevations to pre-covers figure}.
\begin{pf}
Let $D$ be the closure of a component of $\hat{C}_j\smallsetminus\bar{C}_j$. Without loss identify $D\equiv [-1,1]$. Consider the canonical completion $\tilde{X}\to X$ of $\bar{X}\to X$.  There is a unique lift $\tilde{\delta}^+_j:D\to\tilde{X}$ of $\delta_i$ to $\tilde{X}$ so that $\tilde{\delta}^+_j(1)=\bar{\delta}_j(1)$.  Let $\epsilon>0$ be minimal such that $\tilde{\delta}^+_j(\epsilon-1)$ lies in a vertex space. Let $X'\subset\tilde{X}$ be the pre-cover consisting of $\bar{X}$ together with the vertex spaces and edge spaces containing $\tilde{\delta}^+_j((\epsilon-1,1))$.  Let $\tilde{e}_+$ be the edge space of $\tilde{X}$ so that $\tilde{\delta}^+_j([-1,\epsilon-1])\subset \tilde{e}_+\times[-1,1]$ and without loss of generality assume that $\tilde{\delta}^+_j(\epsilon-1)\in\tilde{e}_+\times\{+1\}$.

Similarly, there is a unique lift $\tilde{\delta}^-_j:D\to\tilde{X}$ of $\delta_i$ so that $\tilde{\delta}^-_j(-1)=\bar{\delta}_j(-1)$. Then $\tilde{\delta}^-_j([-1,\epsilon-1])\subset \tilde{e}_-\times[-1,1]$ for some unique edge space $\tilde{e}_-$ of $\tilde{X}$. Without loss of generality, assume $\tilde{\delta}^-_j(-1)\in\tilde{e}_-\times\{-1\}$.

The edge spaces $\tilde{e}_+$ and $\tilde{e}_-$ are both simply connected covers of some edge space $e$ of $X$.  There exists some unique covering transformation $\tau:\tilde{e}_+\to\tilde{e}_-$ such that, whenever $x\in(\epsilon-1,-1)$ with $\tilde{\delta}^+_j(x)\in\tilde{e}_+\times\{\frac{1}{2}\}$,
$$
\tau\circ\tilde{\delta}^+_j(x)=\tilde{\delta}^-_j(x).
$$
Now let $\hat{X}$ be the pre-cover given by $X'$ together with the additional edge space $\tilde{e}_+$; the additional edge maps are $\partial^{\tilde{e}_+}_+$ and $\partial^{\tilde{e}_-}_-\circ\tau$.

Since $\{\bar{\delta}_j\}$ is disparate, this can be done for every of component of $\hat{C}_j\smallsetminus \bar{C}_j$.
\end{pf}

\subsection{Ensuring disparity}

In the light of Lemma \ref{Making elevations full}, we will need to be able to guarantee that collections of elevations are disparate.  We prove under certain assumptions that, in the underlying graph of the canonical completion, any pair of elevations shares at most finitely many vertices.  Therefore, by enlarging our pre-cover to a larger core of the canonical completion, a finite set of elevations can be made disparate.

\begin{lem}[Making elevations disparate]\label{Ensuring disparity}
Suppose $X$ is an $\mathrm{ICE}$ space. Let $\{\delta_i:C_i\to X\}$ be a finite, independent collection of hyperbolic loops such that each $\delta_i$ generates a maximal abelian subgroup of $\pi_1(X)$. Let $\bar{X}$ be a connected pre-cover with finitely generated fundamental group and $\{\bar{\delta}_j:\bar{C}_j\to \bar{X}\}$ a finite collection of elevations of the $\delta_i$. Then there exists a pre-covering $\hat{X}\to X$ extending $\bar{X}\to X$ such that:
\begin{enumerate}
\item $\hat{X}$ has finite underlying graph;
\item $\pi_1(\bar{X})=\pi_1(\hat{X})$;
\item each $\bar{\delta}_j$ extends to an elevation $\hat{\delta}_j$
to $\hat{X}$; and
\item the collection $\{\hat{\delta}_j\}$ is disparate.
\end{enumerate}
\end{lem}
\begin{pf}
Let $H=\pi_1(\bar{X})$, and let $T$ be the Bass--Serre tree of the induced splitting of $\pi_1(X)$, so $\tilde{\Gamma}=H\backslash T$ is the underlying graph of the canonical completion $\tilde{X}$ of $\bar{X}$.  Think of $\tilde{\Gamma}$ as a graph of groups.  The underlying graph $\bar{\Gamma}$ of $\bar{X}$ is a finite core for $\tilde{\Gamma}$. The Bass--Serre tree of $\bar{\Gamma}$ is naturally identified with an $H$-invariant subtree $T^H$ of $T$.

The first observation is that, outside of a finite subgraph, all the edge and vertex stabilizers of $\tilde{\Gamma}$ are trivial. Otherwise, there are elements of $\pi_1(X)$ fixing arbitrarily large subtrees of $T$, contradicting acylindricality.  Therefore, enlarging $\bar{X}$ in $\tilde{X}$, it can be assumed that every edge stabilizer in $\tilde{\Gamma}\smallsetminus\bar{\Gamma}$ is trivial. As a consequence, $H$ freely permutes the components of $T\smallsetminus T^H$.

Each $\bar{\delta}_j$ is the restriction of some elevation $\tilde{\delta}_j:\tilde{C}_j\to\tilde{X}$ which corresponds to a double coset $Hg_j\langle\delta_i\rangle$. Consider the natural map $\phi_{\tilde{X}}:\tilde{X}\to\tilde{\Gamma}$ and set $\lambda_j=\phi_{\tilde{X}}\circ\tilde{\delta}_j$.  By Lemma \ref{Axis in elevation}, the image of $\lambda_j$ contains the projection of the line $l_j=g_j\Axis(\delta_i)\subset T$ to $\tilde{\Gamma}$.

Since $\lambda_j$ is proper (by Lemma \ref{Proper loops}), each $\bar{C}_j$ is compact whenever $\bar{\Gamma}$ is a finite graph.

Note that if $\tilde{\delta}_j$ is infinite-degree then $l_j\cap T^H$  is finite, so there exist two distinct components $\Lambda^\pm_j$ of $T\smallsetminus T^H$ with $l_j\cap\Lambda^\pm_j$ an infinite ray.

We will now show that no two $\lambda_j$ share more than finitely many edges in their images.  Suppose that the intersection of the images of $\lambda_1$ and $\lambda_2$ contains infinitely many edges.  Then $\tilde{\delta}_1$ and $\tilde{\delta}_2$ are infinite-degree elevations and the images of $l_1$ and $l_2$ in $\tilde{\Gamma}$ share infinitely many edges, since $H$ is finitely generated.  Without loss of generality, there exist infinitely many $h_i\in H$ mapping an edge in $\Lambda^+_1\cap l_1$ to an edge in $\Lambda^+_2\cap l_2$.  Since $H$ freely permutes the components of $T\smallsetminus T^H$ these $h_i$ are all equal to some $h\in H$, and $hl_1\cap l_2$ contains infinitely many edges. By acylindricality $\delta_1^{{hg_1}^{-1}}$ and $\delta_2^{g_2^{-1}}$ commute, so $\delta_1=\delta_2$ and
$$
Hg_1\langle \delta_1\rangle=Hg_2\langle\delta_2\rangle.
$$
So $\tilde{\delta}_1$ and $\tilde{\delta}_2$ are isomorphic elevations.

Finally, we also need to check that, for each infinite-degree $\tilde{\delta}_j$, the limits $\lim_{t\to\infty}\lambda_j$ and $\lim_{t\to-\infty}\lambda_j$ differ---otherwise, for any finite core of $\tilde{\Gamma}$, the two endpoints of the restriction of $\tilde{\delta}_j$ would always lie in the same vertex space.   However, if the limits $\lim_{t\to\infty}\lambda_j$ and $\lim_{t\to-\infty}\lambda_j$ coincide then there exists $h\in H$ mapping $l_j\cap\Lambda^+_j$ to $l_j\cap\Lambda^-_j$, so by acylindricality $h\in\langle\delta_i^{{g_j}^{-1}}\rangle$.  This contradicts the assumption that the elevations $\tilde{\delta}_j$ are of infinite degree.

In summary, $\bar{X}$ can be enlarged to a larger core $\hat{X}\subset\tilde{\Gamma}$ satisfying the conclusions of the lemma.
\end{pf}

\section{Constructing finite-sheeted covers}\label{Constructing finite-sheeted covers section}

\subsection{Tame covers}

One might hope to prove Theorems \ref{Theorem A} and \ref{Theorem B} for all groups in $\ICE$, and hence for all limit groups, by induction.  However, neither subgroup separability nor local retractions is strong enough to serve as an inductive hypothesis.  We need a still stronger property---\emph{tameness}.

\begin{defn}\label{Tame covers definition}
Consider a complex $X$, a covering $X'\to X$ and a finite (possibly empty) collection of independent, essential loops
\[
\curlyL=\{\delta_i:C_i\to X\}.
\]
The cover $X'$ is \emph{tame over $\curlyL$} if the following holds.

Let $\Delta\subset X'$ be a finite subcomplex and let
$$
\{\delta'_j:C'_j\to X'\}
$$
be a finite collection of (pairwise non-isomorphic) infinite-degree elevations where each $\delta'_j$ is an elevation of some $\delta_i\in\curlyL$. Then for all sufficiently large positive integers $d$ there exists an intermediate finite-sheeted covering
$$
X'\to\hat{X}\to X
$$
so that:
\begin{enumerate}
\item every $\delta'_j$ descends to an elevation $\hat{\delta}_j$ of degree $d$;
\item the $\hat{\delta}_j$ are pairwise non-isomorphic;
\item \label{Tameness implies separability}$\Delta$ embeds into $\hat{X}$;
\item \label{Tameness implies virtual retract}there is a retraction $\rho:\pi_1(\hat{X})\to\pi_1(X')$ such that
$$
\rho(\hat{\delta}_{j*}(\pi_1(\hat{C}_j)))\subset \delta'_{j*}(\pi_1(C'_j))
$$ for each $j$.
\end{enumerate}
\end{defn}

Suppose that a covering $X'\to X$ is tame over the empty set of loops in $X$.  It follows, by Lemma \ref{Topological interpretation},  from property \ref{Tameness implies separability} of Definition \ref{Tame covers definition} that $\pi_1(X')$ is separable in $\pi_1(X)$.  Likewise, it follows immediately from property \ref{Tameness implies virtual retract} of Definition \ref{Tame covers definition} that $\pi_1(X')$ is a virtual retract of $\pi_1(X)$.  We therefore have the following.

\begin{lem}\label{Tame implies ss and vr}
If a covering of a graph of spaces $X'\to X$ is tame over the empty set of loops in $X$ then $\pi_1(X')$ is separable and a virtual retract in $\pi_1(X)$.
\end{lem}

As a first example, covers of tori are tame over coordinate circles.

\begin{lem}\label{Torus case}
If $T$ is an $n$-torus, $T'\to T$ is a covering and $\delta:S^1\to T$ is a coordinate circle then $T'$ is tame over the singleton $\{\delta\}$.
\end{lem}
\begin{pf}
Identify $T$ with $(\R/\Z)^n$ and let $\delta$ be the loop
$$
t\mapsto (t,0,\ldots,0).
$$
The cover $T'$ can be taken to be the cylinder
$$
\R\times\R^{k-1}\times \R/\mu_{k+1}\Z\times\ldots\times\R/\mu_n\Z
$$
for some integers $\mu_i$, and the elevations $\delta'$ can be taken to be the maps
$$
t\mapsto (t,x_i)
$$
for distinct $x_i\in \R^{k-1}\times\R/\mu_{k+1}\Z\times\ldots\times\R/\mu_n\Z$.  Consider a compact cylindrical subspace of the form
$$
\Omega = [a_1,b_1]\times\ldots\times[a_k,b_k]\times \R/\mu_{k+1}\Z\times\ldots\times\R/\mu_n\Z
$$
for integers $a_i,b_i$.  If the $a_i$ are sufficiently negative and the $b_i$ are sufficiently positive then $\Delta$ lies in the interior of $\Omega$ and, furthermore, $\Omega$ intersects each $\delta_i$ in a non-trivial interval $[a_1,b_1]\times\{x_i\}$.  As long as $d$ is sufficiently large we can take $b_1=a_1+d$. There is an action of $\Z^k$ on $T'$ by covering transformations for which $\Omega$ is a fundamental domain; let $\hat{T}$ be the quotient. This is the required intermediate covering.
\end{pf}

However, if $\curlyL$ consists of two or more loops on $T$ then the covering $T'\to T$ need not be tame over $\curlyL$.  This is what goes wrong in the example given in \cite{BKS} of a 3-manifold group that is not subgroup separable.

\subsection{Making pre-covers finite-sheeted}

Tameness enables us to pass from pre-covers with finite underlying graph to finite-sheeted pre-covers.

\begin{prop}[Passing to finite-sheeted pre-covers]\label{Making pre-covers finite}
Let $X$ be a graph of spaces in which every edge space is a circle, and let $\{\delta_i:C_i\to X\}$ be a collection of elliptic loops confined to the vertex spaces of $X$.  Let $X'\to X$ be a pre-covering with finite underlying graph. Every vertex space $V'$ of $X'$ covers some vertex space $V$ of $X$; assume that each $V'$ is tame over the union of the set of edge maps incident at $V$ and the set of $\delta_i$ whose images lie in $V$. Let $\Delta\subset X'$ be a finite subcomplex and $\{\delta'_j:C'_j\to X'\}$ a collection of infinite-degree elevations of the $\delta_i$. Then for all sufficiently large $d$ there exists a finite-sheeted intermediate pre-covering
$$
X'\to\bar{X}\to X
$$
such that:
\begin{enumerate}
\item the $\{\delta'_j\}$ descend to distinct, full elevations $\{\bar{\delta}_j:\bar{C}_j\to \bar{X}\}$ so that, for each $j$, $\bar{C}_j\to C$ is a covering of degree $d$;
\item $\Delta$ embeds into $\bar{X}$; and
\item there exists a retraction $\rho:\pi_1(\bar{X})\to\pi_1(X')$ such that
$$
\rho(\bar{\delta}_{j*}(\pi_1(\bar{C}_j)))\subset \delta'_{j*}(\pi_1(C'_j))
$$
for each $j$.
\end{enumerate}
\end{prop}
\begin{pf}
Enlarging $\Delta$ it can be assumed that it contains every compact edge space of $X'$.  Let $V'$ be a vertex space of $X'$ covering the vertex $V$ of $X$.  Set $\Delta_{V'}=V'\cap\Delta$. Consider the edge maps $\partial'_i:e_i\to V'$ of $X'$ incident at $V'$ that are infinite-degree elevations of $\partial_\pm:e\to V$. By the hypothesis of tameness there exists an intermediate finite-sheeted covering
$$
V'\to\bar{V}\to V
$$
so that $\Delta_{V'}$ embeds in $\bar{V}$, each $\partial'_i$ descends to a degree $d$ elevation
$$
\bar{\partial}_i:\bar{e}\to \bar{V}
$$
of $\partial_\pm$ and each $\delta'_j$ with image in $V'$ descends to a degree $d$ elevation $\bar{\delta}_j:\bar{C}_j\to\bar{V}$. If $d$ is sufficiently large then it can be taken to be the same $d$ over all vertex spaces of $X'$.

Let $\bar{X}$ be the graph of spaces with the same underlying graph as $X'$ but in which each vertex space $V'$ is replaced by the corresponding $\bar{V}$.  If $e'$ is an edge space of $X'$ then the edge map
$$
\partial'_\pm:e'\to V'
$$
descends to a finite-degree map $\bar{\partial}_\pm:\bar{e}_\pm\to \bar{V}$.  By construction $\bar{e}_+$ covers $e$ with the same degree as $\bar{e}_-$ so this information defines a finite-sheeted pre-cover $\bar{X}$ as required.

It is clear that $\Delta$ embeds into $\bar{X}$.  As for the retraction, $\pi_1(\bar{X})$ decomposes as a graph of groups, with the same underlying graph as the decomposition for $\pi_1(X')$, and every vertex group of $\pi_1(X')$ is a retract of the corresponding vertex group of $\pi_1(\bar{X})$.  These retractions piece together to give the required retraction $\rho:\pi_1(\bar{X})\to\pi_1(X')$.
\end{pf}

\subsection{Completing pre-covers to covers}

Consider an $\ICE$ space $X$, which as usual is constructed by gluing a torus $T$ to a simpler $\ICE$ space $Y$.  Let $H\subset\pi_1(X)$ be a finitely generate subgroup and $X^H\to X$ the corresponding covering.  Proposition \ref{Making pre-covers finite} will enable us, under suitable hypotheses, to pass from a core of $X'$ to a finite-sheeted pre-cover $\bar{X}$. The next step is to add vertex spaces so that we can apply Stallings' Principle. Furthermore, we would like to do this in such a way that the resulting cover retracts onto the pre-cover $\bar{X}$. To this end, we will need to construct finite-sheeted covers of the vertex spaces with special properties. The next lemma does this for $T$.

\begin{lem}\label{Controlling covers of tori}
Let $T$ be a torus and $\delta:S^1\to T$ an essential loop. Then for every positive integer $d$ there exists a finite-sheeted covering $\hat{T}_d\to T$ so that $\delta$ has a single elevation $\hat{\delta}$ to $\hat{T}_d$, and $\hat{\delta}$ is of degree $d$.
\end{lem}
\begin{pf}
Since $\pi_1(T)$ is abelian we can assume that $\delta$ is a based loop.  Let $\pi_1(S^1)=\langle t\rangle$ and consider the decomposition
$$
\pi_1(T)=C\oplus A
$$
where $C$ is maximal cyclic with $\delta_*t\in C$ and $A$ is free abelian. Setting $\pi_1(\hat{T}_d)=\langle d\delta_*t\rangle\oplus A$ gives the required covering.
\end{pf}
Note that $\pi_1(\hat{T}_d)$ retracts onto $\langle\hat{\delta}\rangle=\langle d\delta_*t\rangle$. For $Y$, the inductive hypothesis that $\pi_1(Y)$ admits local retractions will ensure the existence of the necessary covers.

\begin{lem}\label{Controlling LR covers}
Let $Y$ be a space such that $\pi_1(Y)$ has local retractions and let $\delta:S^1\to Y$ be a based essential closed loop.  Then for every positive integer $d$ there exists a finite-sheeted covering $\hat{Y}_d\to Y$ so that $\delta$ has an elevation $\hat{\delta}$ to $\hat{Y}_d$ of degree $d$ and $\pi_1(\hat{Y}_d)$ retracts onto $\langle\hat{\delta}\rangle$.
\end{lem}
\begin{pf}
Since $\pi_1(Y)$ admits local retractions, there exists a finite-sheeted covering $\hat{Y}_d\to Y$ so that $\pi_1(\hat{Y}_d)$ retracts onto $\langle\delta^d\rangle$.  Note that $\delta^k\notin\pi_1(\hat{Y}_d)$ whenever $0<k<d$.  Hence the double coset $\pi_1(\hat{Y}_d)1\langle\delta\rangle$ corresponds to an elevation of $\delta$ of degree $d$.
\end{pf}

If we don't require $\pi_1(\hat{Y}_d)$ to retract onto $\langle\hat{\delta}\rangle$ then Lemma \ref{Controlling LR covers} can be significantly strengthened for residually free groups. Indeed, whenever the fundamental group of a space $Y$ is residually torsion-free nilpotent, the cover $\hat{Y}_d$ can be taken to be normal, so every elevation of $\delta$ to $\hat{Y}_d$ is of degree $d$  (see Remark 3.1 of \cite{Gi1}).  By a theorem of W.~Magnus, free groups (and hence residually free groups) are residually torsion-free nilpotent. However, we will not need this strengthening.

\begin{prop}[Completing a finite-sheeted pre-cover to a cover]\label{Completing pre-covers to covers}
Let $X$ be an $\ICE$ space constructed by gluing together a torus $T$ and a simpler $\ICE$ space $Y$, as above.  Assume that $\pi_1(Y)$ admits local retractions. Let $\bar{X}\to X$ be a finite-sheeted pre-covering.  Then there exists an inclusion $\bar{X}\hookrightarrow\hat{X}$ extending $\bar{X}\to X$ to a covering $\hat{X}\to X$ such that $\pi_1(\hat{X})$ retracts onto $\pi_1(\bar{X})$.
\end{prop}
\begin{pf}
Let $d$ be a positive integer.  By Lemma \ref{Controlling LR covers} there exists a finite-sheeted covering $\bar{Y}_d\to Y$ with an elevation of $\partial_+$ of degree $d$.  Likewise, by Lemma \ref{Controlling covers of tori} there exists a finite-sheeted covering $\bar{T}_d\to T$ with just one elevation of $\partial_-$, and that elevation is of degree $d$.

Let $m_d$ be the number of hanging elevations of $\partial_-$ to $\bar{X}$ of degree $d$. Consider the pre-cover
$$
\bar{X}_+=\bar{X}\sqcup\coprod_{d}\coprod_{i=1}^{m_d} \bar{Y}_d.
$$
Let $p_d$ be the number of degree $d$ elevations of $\partial_+$ to $\bar{X}_+$ and $q_d$ the number of degree $d$ elevations of $\partial_-$ to $\bar{X}^+$.  By construction, $n_d=p_d-q_d\geq 0$.  Pairing up every degree $d$ hanging elevation of $\partial_-$ in $\bar{X}$ with a hanging elevation of $\partial_+$ in some $\bar{Y}_d$, we can assume that $\bar{X}_+$ is connected and has $n_d$ hanging elevations of $\partial_+$ of degree $d$. Consider the pre-cover
$$
\bar{X}_0=\bar{X}_+\sqcup\coprod_{d}\coprod_{i=1}^{n_d} \bar{T}_d
$$
By construction, there exist $n_d$ degree $d$ hanging elevations of $\partial_+$ to $\bar{X}_0$ and $n_d$ degree $d$ elevations of $\partial_-$ to $\bar{X}_0$.  So, by Stallings' Principle, $\bar{X}_0$ extends to a finite-sheeted cover $\hat{X}$ of $X$.

Each $\pi_1(\bar{T}_d)$ retracts onto the edge group to which it was glued; then each $\pi_1(\bar{Y}_d)$ retracts onto the edge group to which it was glued.  So $\pi_1(\hat{X})$ retracts onto $\pi_1(\bar{X})$.
\end{pf}

\subsection{The induction}

In this section we prove our main theorem.

\begin{thm}\label{Main theorem}
Let $X$ be an $\ICE$ space, let $H\subset\pi_1(X)$ be a finitely generated subgroup and let $X^H\to X$ be the corresponding covering.  Suppose $\curlyL$ is a (possibly empty) finite, independent set of loops that generate maximal abelian subgroups of $\pi_1(X)$. Then $X^H$ is tame over $\curlyL$.
\end{thm}
The proof of Theorem \ref{Main theorem} is an induction.  The inductive step is similar to the proof of Corollary \ref{Motivating example}, which provides the base step.

\begin{pfof}{Theorem \ref{Main theorem}}
Fix a finitely generated subgroup $H\subset\pi_1(X)$ and let $X^H\to X$ be the corresponding covering space.  For notational convenience, divide $\curlyL$ into hyperbolic loops $\{\delta_i:D_i\to X\}$ and elliptic loops $\{\epsilon_i:E_i\to X\}$.   After modifying by a homotopy, we can assume that each $\epsilon_i$ is contained in a vertex of space of $X$.  Let $\Delta\subset X^H$ be a finite subcomplex, and as in the definition of tameness let $\{\delta^H_j:D^H_j\to X^H\}$ be a finite collection of infinite-degree elevations of the $\{\delta_i\}$ and $\{\epsilon'_j:E'_j\to X^H\}$ a finite collection of infinite-degree elevations of the $\{\epsilon_i\}$.

Since $H$ is finitely generated, there exists a core $X'\subset X^H$. Enlarging $X'$ if necessary, we can assume that $\Delta\subset X'$, that each $\delta^H_j$ restricts to an elevation $\delta'_j:D'_j\to X'$ and that the image of each $\epsilon'_j$ is contained in $X'$.

By Lemma \ref{Ensuring disparity} we can enlarge $X'$ further and assume that the $\delta'_j$ are disparate. For all sufficiently large $d$, since $D'_j$ is a finite collection of compact intervals, there exists an embedding into a circle $D'_j\to \bar{D}_j$ so that the $d$-sheeted covering $\bar{D}_j\to D_i$ extends $D'_j\to D_i$; so by Lemma \ref{Making elevations full} we can extend $X'$ to a pre-cover $\bar{X}$ such that each $\delta'_j$ extends to a full elevation
$$
\bar{\delta}_j:\bar{D}_j\to\bar{X}.
$$
Enlarging $\Delta$ again, we can assume that the images of the $\bar{\delta}_j$ are contained in $\Delta$.  By Proposition \ref{Making pre-covers finite} there exists some intermediate finite-sheeted pre-covering
$$
\bar{X}\to\hat{X}\to X
$$
into which $\Delta$ injects, and such that each $\epsilon'_j$ descends to a full elevation $\hat{\epsilon}_j:\hat{E}_j\to\hat{X}$ with $\hat{E}_j\to E_i$ a covering of degree $d$. Since $\Delta$ injects into $\hat{X}$ we already have that $\bar{\delta}_j$ descends to an elevation
$$
\hat{\delta}_j:\bar{D}_j\to\hat{X}.
$$
Finally, $\hat{X}$ can be extended to a finite-sheeted covering $\hat{X}^+$ by Proposition \ref{Completing pre-covers to covers}.

By construction, $\Delta$ injects into $\hat{X}^+$ and $\pi_1(\hat{X}^+)$ retracts onto $\pi_1(\hat{X})$, which in turn retracts onto $\pi_1(\bar{X})$.  By Lemma \ref{Making elevations full}, $H=\pi_1(X')$ is a free factor in $\pi_1(\bar{X})$.  So, in summary, there exists a retraction $\rho:\pi_1(\hat{X}^+)\to H$.  Furthermore, $\rho$ is either the identity  or the trivial homomorphism on the cyclic subgroups generated by each $\hat{\delta}_j$ and $\hat{\epsilon}_j$.
\end{pfof}

\subsection{Conclusion}\label{QI subgroups}

Theorems \ref{Theorem A} and \ref{Theorem B} follow easily.  By Lemma \ref{Tame implies ss and vr} and Theorem \ref{Main theorem}, groups in $\ICE$ are subgroup separable and have local retractions.   Applying Theorem \ref{Limit groups embed in ICE groups} and the facts that subgroup separability and local retractions pass to subgroups, it follows that limit groups enjoy both these properties.

\begin{cor}[Theorem \ref{Theorem A}]\label{LGs are subgroup separable}
Limit groups are subgroup separable.
\end{cor}

\begin{cor}[Theorem \ref{Theorem B}]\label{LGs have local retractions}
Limit groups have local retractions.
\end{cor}

It follows from Corollary \ref{LGs are subgroup separable} that limit groups have solvable generalized word problem.

\begin{cor}
Limit groups have solvable generalized word problem.
\end{cor}

Corollary \ref{LGs have local retractions} and Lemma \ref{LR implies qi embedded subgroups} together imply that finitely generated subgroups of limit groups are quasi-isometrically embedded.

\begin{cor}\label{QI embedded subgroups}
Every finitely generated subgroup of a limit group is quasi-isometrically embedded.
\end{cor}

F.~Dahmani proved in \cite{Da} the stronger result that finitely generated subgroups of limit groups are quasi-convex in the sense of convergence groups. In contrast, our techniques are purely topological, and indeed the notion of quasi-convexity is not well-defined without introducing some geometry.

In this context one can deduce more from Corollary \ref{QI embedded subgroups}.  Emina Alibegovic and Mladen Bestvina showed in \cite{AB04} that any limit group $G$ acts geometrically on a \emph{CAT(0) space with isolated flats} $X$. Christopher Hruska proved that in this situation a finitely generated subgroup $H$ of $G$ is quasi-convex with respect to the action on $X$ if and only if $H$ is quasi-isometrically embedded in $G$ (Theorem 1.1 of \cite{Hr}, see also Proposition 4.1.6 of \cite{HK}). Combining these results with Corollary \ref{QI embedded subgroups} we obtain a quasi-convexity result.

\begin{cor}
Let $G$ be a limit group with an action $\rho$ by isometries on a CAT(0) space with isolated flats.  Then every finitely generated subgroup $H$ of $G$ is quasi-convex with respect to the action $\rho$.
\end{cor}

The proof of Theorem \ref{Main theorem} has deeper consequences for the separability properties of $\ICE$ groups.

\begin{cor}
Let $G\in\ICE$ and let $H$ be a vertex or edge group of $G$.  Then $G$ is coset separable with respect to $H$.
\end{cor}
\begin{pf}
As usual, think of $G$ as the fundamental group of an $\ICE$ space $X$.  We will prove the corollary when $H$ is a vertex group of $G$, corresponding to a vertex space $V$ of $X$.  The proof when $H$ corresponds to an edge group is the same.

Let $H'\subset G$ be a finitely generated subgroup and let $X^{H'}$ be the corresponding cover of $X$.  Assume that $HgH'\neq HhH'$.  Translating this to the language of elevations, this is equivalent to asserting that the elevations of the inclusion map $V\hookrightarrow X$ to the cover $X^{H'}$ are non-isomorphic (in particular, they have disjoint images).  The image of each of these elevations is precisely a vertex space of $X^{H'}$.  Let $V'_g$ be the vertex space corresponding to $HgH'$ and let $V'_h$ be the vertex space corresponding to $HhH'$.

Let $X'$ be a core for $X^{H'}$.  Enlarging $X'$ if necessary, we can assume that $V'_g$ and $V'_h$ are both contained in $X'$.  The proof of Theorem \ref{Main theorem} constructs an intermediate finite-sheeted cover $\hat{X}\to X$ in which distinct vertex and edge spaces of $X'$ descend to distinct vertex and edge spaces of $\hat{X}$.  In particular, $Hg\pi_1(\hat{X})\neq Hh\pi_1(\hat{X})$, as required.
\end{pf}

\bibliographystyle{plain}

\begin{thebibliography}{10}

\bibitem{AB04}
Emina Alibegovi{\'c} and Mladen Bestvina.
\newblock Limit groups are {$\rm CAT(0)$}.
\newblock {\em J. London Math. Soc. (2)}, 74(1):259--272, 2006.

\bibitem{BrHo}
Martin~R. Bridson and James Howie.
\newblock Normalisers in limit groups.
\newblock {\em Math. Ann.}, 337(2):385--394, 2007.

\bibitem{BKS}
R.~G. Burns, A.~Karrass, and D.~Solitar.
\newblock A note on groups with separable finitely generated subgroups.
\newblock {\em Bull. Austral. Math. Soc.}, 36(1):153--160, 1987.

\bibitem{CG04}
Christophe Champetier and Vincent Guirardel.
\newblock Limit groups as limits of free groups.
\newblock {\em Israel J. Math.}, 146:1--75, 2005.

\bibitem{Da}
Fran{\c{c}}ois Dahmani.
\newblock Combination of convergence groups.
\newblock {\em Geom. Topol.}, 7:933--963 (electronic), 2003.

\bibitem{Gi1}
Rita Gitik.
\newblock Graphs and separability properties of groups.
\newblock {\em J. Algebra}, 188(1):125--143, 1997.

\bibitem{Hall49}
Marshall Hall, Jr.
\newblock Subgroups of finite index in free groups.
\newblock {\em Canadian J. Math.}, 1:187--190, 1949.

\bibitem{Ha02}
Allen Hatcher.
\newblock {\em Algebraic Topology}.
\newblock Cambridge University Press, 2002.

\bibitem{Hr}
G.~C. Hruska.
\newblock Geometric invariants of spaces with isolated flats.
\newblock {\em Topology}, 44(2):441--458, 2005.

\bibitem{HK}
G.~Christopher Hruska and Bruce Kleiner.
\newblock Hadamard spaces with isolated flats.
\newblock {\em Geom. Topol.}, 9:1501--1538 (electronic), 2005.
\newblock With an appendix by the authors and Mohamad Hindawi.

\bibitem{KM98a}
O.~Kharlampovich and A.~Miasnikov.
\newblock Irreducible affine varieties over a free group. {I}. {I}rreducibility
  of quadratic equations and {N}ullstellensatz.
\newblock {\em J. Algebra}, 200(2):472--516, 1998.

\bibitem{KM98b}
O.~Kharlampovich and A.~Miasnikov.
\newblock Irreducible affine varieties over a free group. {II}. {S}ystems in
  triangular quasi-quadratic form and description of residually free groups.
\newblock {\em J. Algebra}, 200(2):517--570, 1998.

\bibitem{KM4}
O.~Kharlampovich and A.~Miasnikov.
\newblock Effective {JSJ} decompositions (paper 4 on the theory of a free
  group).
\newblock In {\em Groups, languages, algorithms}, volume 378 of {\em Contemp.
  Math.}, pages 87--212. Amer. Math. Soc., Providence, RI, 2005.

\bibitem{KM3}
O.~Kharlampovich and A.~Miasnikov.
\newblock Implicit function theorem over free groups (paper 3 on the theory of
  a free group).
\newblock {\em J. Algebra}, 290(1):1--203, 2005.

\bibitem{KM5}
O.~Kharlampovich and A.~Miasnikov.
\newblock Elementary theory of free nonabelian groups (paper 5 on the theory of
  a free group).
\newblock {\em J. Algebra}, 302(2):451--552, 2006.

\bibitem{LR2}
D.~D. Long and A.~W. Reid.
\newblock Subgroup separability and virtual retractions of groups.
\newblock \emph{Topology}, to appear, 2006.

\bibitem{Rem89}
V.~N. Remeslennikov.
\newblock {$\exists$}-free groups.
\newblock {\em Sibirsk. Mat. Zh.}, 30(6):193--197, 1989.

\bibitem{Sc1}
Peter Scott.
\newblock Subgroups of surface groups are almost geometric.
\newblock {\em J. London Math. Soc. (2)}, 17(3):555--565, 1978.

\bibitem{SW}
Peter Scott and Terry Wall.
\newblock Topological methods in group theory.
\newblock In {\em Homological group theory (Proc. Sympos., Durham, 1977)},
  volume~36 of {\em London Math. Soc. Lecture Note Ser.}, pages 137--203.
  Cambridge Univ. Press, Cambridge, 1979.

\bibitem{Se1}
Z.~Sela.
\newblock Diophantine geometry over groups. {I}. {M}akanin-{R}azborov diagrams.
\newblock {\em Publ. Math. Inst. Hautes \'Etudes Sci.}, (93):31--105, 2001.

\bibitem{Se2}
Z.~Sela.
\newblock Diophantine geometry over groups. {II}. {C}ompletions, closures and
  formal solutions.
\newblock {\em Israel J. Math.}, 134:173--254, 2003.

\bibitem{Se4}
Z.~Sela.
\newblock Diophantine geometry over groups. {IV}. {A}n iterative procedure for
  validation of a sentence.
\newblock {\em Israel J. Math.}, 143:1--130, 2004.

\bibitem{Se3}
Z.~Sela.
\newblock Diophantine geometry over groups. {III}. {R}igid and solid solutions.
\newblock {\em Israel J. Math.}, 147:1--73, 2005.

\bibitem{Se5a}
Z.~Sela.
\newblock Diophantine geometry over groups. {$\rm V\sb 1$}. {Q}uantifier
  elimination. {I}.
\newblock {\em Israel J. Math.}, 150:1--197, 2005.

\bibitem{Se6}
Z.~Sela.
\newblock Diophantine geometry over groups. {VI}. {T}he elementary theory of a
  free group.
\newblock {\em Geom. Funct. Anal.}, 16(3):707--730, 2006.

\bibitem{Se5b}
Z.~Sela.
\newblock Diophantine geometry over groups {$V\sb 2$}: quantifier elimination.
  {II}.
\newblock {\em Geom. Funct. Anal.}, 16(3):537--706, 2006.

\bibitem{SeQ}
Zlil Sela.
\newblock Diophantine geometry over groups: a list of research problems.
\newblock \texttt{http://www.ma.huji.ac.il/\textasciitilde zlil/problems.dvi}.

\bibitem{S77}
Jean-Pierre Serre.
\newblock {\em Arbres, amalgames, {${\rm SL}\sb{2}$}}.
\newblock Soci\'et\'e Math\'ematique de France, Paris, 1977.
\newblock Avec un sommaire anglais, R\'edig\'e avec la collaboration de Hyman
  Bass, Ast\'erisque, No. 46.

\bibitem{St}
John~R. Stallings.
\newblock Topology of finite graphs.
\newblock {\em Invent. Math.}, 71(3):551--565, 1983.

\bibitem{Wilt1}
Henry Wilton.
\newblock Elementarily free groups are subgroup separable.
\newblock To appear in \emph{Proc. London Math. Soc.} ArXiv:~math.GR/0511414.

\bibitem{Wilt2}
Henry Wilton.
\newblock Solutions to {B}estvina and {F}eighn's exercises on limit groups.
\newblock To appear in \emph{Geometry and Cohomology in Group Theory}
  \emph{(Durham, 2003)}. ArXiv: math.GR/0604137.

\bibitem{Wilt3}
Henry Wilton.
\newblock Subgroup separability of limit groups.
\newblock \emph{PhD thesis, Univ. London}, 2006.

\bibitem{Wise1}
Daniel~T. Wise.
\newblock Subgroup separability of graphs of free groups with cyclic edge
  groups.
\newblock {\em Q. J. Math.}, 51(1):107--129, 2000.

\end{thebibliography}

\bigskip\bigskip\centerline{\textbf{Author's address}}

\smallskip\begin{center}\begin{tabular}{l}%
Department of Mathematics\\
1 University Station C1200\\
Austin, TX 78712-0257\\
USA\\
{\texttt{henry.wilton@math.utexas.edu}}
\end{tabular}\end{center}

\end{document}